\documentclass[MSNbibl,number,citesort,dvips]{arxbj}

\aid{0}
\volume{20}
\issue{2}
\pubyear{2014}
\firstpage{416}
\lastpage{456}
\doi{10.3150/12-BEJ492} 

\makeatletter
\newproclaim{assumption}{Assumption}
\newproclaim{definition}{Definition}
\newcommand{\rrVert}{\Vert}
\newcommand{\rrvert}{\vert}
\newcommand{\llVert}{\Vert}
\newcommand{\llvert}{\vert}
\newcommand{\eqref}[1]{(\ref{#1})}
\newtheorem{theorem}{Theorem}
\newtheorem{lemma}{Lemma}
\newtheorem{proposition}{Proposition}
\newremark{remark}{Remark}
\makeatother

\begin{document}
\begin{frontmatter}

\title{Total variation approximations and conditional limit theorems
for multivariate
regularly varying random walks conditioned~on~ruin}
\runtitle{Conditional limit theorems for multivariate regularly
varying random walks}

\begin{aug}
\author[1]{\fnms{Jose} \snm{Blanchet}\thanksref{1}}
\and
\author[2]{\fnms{Jingchen} \snm{Liu}\corref{}\thanksref{2}\ead[label=e2]{jcliu@stat.columbia.edu}}
\runauthor{J. Blanchet and J. Liu} 
\address[1]{IEOR Department, Columbia University, 500 West 120th Street,
340 W. Mudd Building,
New York, NY 10027, USA}
\address[2]{Department of Statistics, Columbia University, 1255 Amsterdam Ave, New York, NY 10027,
USA.
\printead{e2}}
\end{aug}

\received{\smonth{9} \syear{2011}}
\revised{\smonth{7} \syear{2012}}

%
\begin{abstract}
We study a new technique for the asymptotic analysis of heavy-tailed systems
conditioned on large deviations events. We illustrate our approach in the
context of ruin events of multidimensional regularly varying random
walks. Our
approach is to study the Markov process described by the random walk
conditioned on hitting a rare target set. We construct a Markov chain whose
transition kernel can be evaluated directly from the increment
distribution of
the associated random walk. This process is shown to approximate the
conditional process of interest in total variation. Then, by analyzing the
approximating process, we are able to obtain asymptotic conditional joint
distributions and a conditional functional central limit theorem of several
objects such as the time until ruin, the whole random walk prior to
ruin, and
the overshoot on the target set. These types of joint conditional limit
theorems have been obtained previously in the literature only in the one
dimensional case. In addition to using different techniques, our results
include features that are qualitatively different from the one dimensional
case. For instance, the asymptotic conditional law of the time to ruin
is no
longer purely Pareto as in the multidimensional case.
\end{abstract}


\end{frontmatter}

\section{Introduction}\label{sec1}

The focus of this paper is the development of a precise asymptotic description
of the distribution of a multidimensional regularly varying random walk (the
precise meaning of which is given in Section \ref{sectionIS})
conditioned on
hitting a rare target set represented as the union of half spaces. In
particular, we develop tractable total variation approximations (in the sample
path space), based on change-of-measure techniques, for such conditional
stochastic processes. Using these approximations we are able to obtain,
as a
corollary, joint conditional limit theorems of specific objects such as the
time until ruin, a Brownian approximation up to (just before) the time of
ruin, the ``overshoot'', and the
``undershoot''. This is the first
paper, as
far as we know, that develops refined conditional limit theorems in a
multidimensional ruin setting; results in one dimensional settings include,
for instance, \cite{AsmKlup96,KlupKypMall04,DonKyp06}; see also \cite
{AsFo12} for extensions concerning regenerative processes.

The techniques developed to obtain our results are also different from those
prevalent in the literature and interesting qualitative features arise in
multidimensional settings. For instance, surprisingly, the asymptotic
conditional time to ruin is no longer purely Pareto, as in the one dimensional
case. A slowly-varying correction is needed in the multidimensional setting.

Other results in the case of multidimensional regularly varying random walks
have been obtained by using a weak convergence approach (see
\cite{HultLinMikSam05}). In contrast to the weak convergence approach, which
has been applied to non-Markovian settings \cite{HultGenna10}, the
techniques that we present here appear to be suited primarily to Markovian
settings. On the other hand, the approach that we shall demonstrate
allows to
obtain finer approximations to conditional objects, for instance, total
variation approximations, and conditional central limit theorems. Moreover,
the present approach has also been applied to non-regularly varying
heavy-tailed settings, at least in one dimension (see \cite
{BL2012SPA}). In
addition, if there is a need to improve upon the quality of the approximations
the method that we advocate readily provides Monte Carlo algorithms
that can
be shown to be optimal in a sense of controlling the relative mean squared
error uniformly in the underlying large deviations parameter (see
\cite{BLMultHT10}). Standard approximation techniques for heavy-tailed large
deviations cannot be easily translated into efficient Monte Carlo algorithms
(see \cite{pASM97a}).

The multidimensional problem that we consider here is a natural
extension of
the classical one dimensional random ruin problem in a so-called
renewal risk
model (cf. the texts of \cite{ASM00Ruin,pASM03bk}). We consider a
$d$-dimensional regularly varying random walk $S= ( S_{n}\dvt n\geq
1 )
$ with $S_{0}=0$ and drift $\eta\in\mathbb{R}^{d}$ so that
$ES_{n}=n\eta\neq
0$. Define
\[
T_{bA^{\ast}}=\inf\bigl\{n\geq0\dvt S_{n}\in bA^{\ast}\bigr
\},
\]
where $A^{\ast}$ is the union of half spaces and $\eta$ points to the interior
of some open cone that does not intersect $A^{\ast}$. The paper
\cite{HultLin06} notes that $P(T_{bA^{\ast}}<\infty)$ corresponds to
the ruin
probabilities for insurance companies with several lines of business. Using
natural budget constraints related to the amount of money that can be
transferred from one business line to another, it turns out that the target
set takes precisely the form of the union of half spaces as we consider here.

Our goal is to illustrate new techniques that can be used to describe very
precisely the conditional distribution of the heavy-tailed processes given
$T_{bA^{*}}<\infty$. Our approximations allow to obtain, in a
relatively easy
way, extensions of previous results in the literature that apply only
in the
one dimensional settings. Asmussen and Kl\"{u}ppelberg \cite{AsmKlup96}
provide conditional limit theorems for the overshoot, the undershoot,
and the
time until ruin given the eventual occurrence of ruin. Similar results have
been obtained recently in the context of Levy processes (see
\cite{KlupKypMall04} and references there in). We apply our results
here to
obtain multidimensional analogues of their conditional limit theorems and
additional refinements, such as conditional central limit theorems.

Our general strategy is based on the use of a suitable change of
measure and
later on coupling arguments. The idea is to approximate the conditional
distribution of the underlying process step-by-step, in a Markovian
way, using
a mixture of a large increment that makes the random walk hit the
target set
and an increment that follows the nominal (original) distribution. The mixture
probability is chosen depending on the current position of the random walk.
Intuitively, given the current position, the selection of the mixture
probability must correspond to the conditional probability of reaching the
target set in the immediate next step given that one will eventually
reach the
target set. The conditional distribution itself is also governed by a Markov
process, so the likelihood ratio (or Radon--Nikodym derivative) between the
conditional distribution and our approximating distribution can be explicitly
written in terms of the ratios of the corresponding Markov kernels. By showing
that the second moment of the likelihood ratio between the conditional
distribution and our approximating distribution approaches unity as the rarity
parameter $b\rightarrow\infty,$ we can reach the desired total variation
approximation (cf. Lemma \ref{TV1}). A crucial portion of our strategy
involves precisely obtaining a good upper bound on the second moment of the
likelihood ratio. The likelihood ratio is obtained out of a Markovian
representation of two measures. It is natural to develop a Lyapunov-type
criterion for the analysis of the second moment. This approach is
pursued in
Section \ref{SecLBTV}, where we introduce our Lyapunov criterion and develop
the construction of the associated Lyapunov function which allows us to bound
the second moment of the likelihood ratio of interest as a function of the
initial state of the process.

There are several interesting methodological aspects of our techniques that
are worth emphasizing. First, the application of change-of-measure
ideas is
common in the light-tailed settings. However, it is not at all standard in
heavy-tailed settings. A second interesting methodological aspect of our
technique is the construction of the associated Lyapunov function. This step
often requires a substantial amount of ingenuity. In the heavy-tailed setting,
as we explain in Section~\ref{SecLBTV}, we can take advantage of the fluid
heuristics and asymptotic approximations for this construction. This approach
was introduced in \cite{BG07AAP} and has been further studied in
\cite{BGL07Questa,BLMultHT10} and \cite{BLldrws07}. In fact,
many of
our ideas are borrowed from the rare-event simulation literature, which
is not
surprising given that our strategy involves precisely the construction
of a
suitable change of measure and these types of constructions, in turn,
lie at
the heart of importance sampling techniques. The particular change of measure
that we use is inspired by the work of \cite{DULEWA06} who applied it
to the
setting of one dimensional (finite) sums of regularly varying increments.

The approximation constructed in this paper is tractable in the sense
that it
is given by a Markovian description which is easy to describe and is explicit
in terms of the increment distribution of the associated random walk (see
(\ref{DKt}) and Theorem \ref{ThmMain1}). This tractability property
has useful
consequences both in terms of the methodological techniques and
practical use.
From a methodological standpoint, given the change-of-measure that we
use to
construct our Markovian description (basically following (\ref
{Dka})), the
result of convergence in total variation provides a very precise justification
of the intuitive mechanism which drives the ruin in the heavy-tailed
situations. In addition, the result allows to directly use this intuitive
mechanism to provide functional probabilistic descriptions that, while less
precise than total variation approximations, emphasize the most important
elements that are present at the temporal scales at which ruin is
expected to
occur (if it happens at all). These functional results are given in Theorem
\ref{ThmMain2}.
Our total variation approximation (Theorem \ref{ThmMain1}) allows to construct
a very natural coupling which makes the functional probabilistic descriptions
given in Theorem \ref{ThmMain2} relatively straightforward in view of standard
strong Brownian approximation results for random walks.

The tractability of our total variation approximation allows for a deep study
of the random walk conditioned on bankruptcy, as mentioned earlier, via
efficient Monte Carlo simulation, at scales that are finer than those provided
by the existing functional limit theorems (cf. \cite
{HultLinMikSam05}). Using
the techniques that we pursue here, the results in \cite{BL2012SPA} establish
\textit{necessary and sufficient} conditions for optimal estimation of
conditional expectations given bankruptcy using importance sampling;
surprisingly, one can ensure finite expected termination time and
asymptotically optimal relative variance control even when the zero variance
change of measure has infinite expected termination time.

As mentioned earlier, we believe that the techniques that we consider
here can
find potential applications in Markovian settings beyond random walks.
This is
a research avenue that we are currently exploring; see, for instance,
\cite{BGL08gg2}, where we apply similar techniques to multi-queues.
Additional results will be reported in the future.

The rest of the paper is organized as follows. In Section \ref
{sectionIS} we
explain our assumptions, describe our approximating process, and state our
main results. The estimates showing total variation approximation,
which are
based on the use of Lyapunov inequalities, are given in Section \ref{SecLBTV}.
Finally, the development of conditional functional central limit
theorems is
given in Section \ref{SectCCLT}.

\section{Problem setup and main results}\label{sectionIS}

\subsection{Problem setup}\label{SecISPS}

Let $(X_{n}\dvt n\geq1)$ be a sequence of independent and identically distributed
(i.i.d.) regularly varying random vectors taking values in $\mathbb{R}^{d}$.
Let $X$ be a generic random variable equal in distribution to $X_{i}$. The
random vector $X$ is said to have a multivariate regularly varying
distribution if there exists a sequence $\{a_n\dvt n \geq1\}$,
$0<a_n\uparrow\infty$, and a non-null Random measure $\mu$ on the
compactified and punctured space $\overline{\mathbb{R}}^{d}\setminus
\{0\}$ with $\mu(\overline{\mathbb{R}}^{d}\setminus{\mathbb
{R}}^{d})=0$ such that, as $n\to\infty$\vspace*{-2.4pt}
%
\begin{equation}
n P\bigl(a_n^{-1} X\in\cdot\bigr) \stackrel v \rightarrow
\mu(\cdot),\label{RVDist}
\end{equation}
where ``$\stackrel v\rightarrow$'' refers to vague convergence.
It can be shown that as $b\to\infty$,\vspace*{-2.4pt}
\[
\frac{P(X\in b \cdot)}{P(\Vert X\Vert_2>b)}\stackrel v \to c \mu (\cdot)
\]
for some $c>0$ (\cite{HultLinMikSam05}, Remark 1.1). To simplify
notation, $a_n$ is chosen such that $nP(\Vert X\Vert_2>a_n)\to1$ and
with this choice of $a_n$ we have that $c=1$.

The random vector $X$ has a relatively very small probability of
jumping into
sets for which $\mu ( B ) =0$. If $P ( \llVert
X\rrVert _{2}>b ) =b^{-\alpha}L ( b ) $ for some
$\alpha>0$ and a slowly varying function $L ( \cdot ) $ (i.e.,
$L ( tb ) /L ( b ) \longrightarrow1$ as
$b\uparrow
\infty$ for each $t>0$), then we say that $\mu ( \cdot )
$ has
(regularly varying) index~$\alpha$. For further information on multivariate
regular variation see \cite{RES06}; the definition provided above corresponds
to the representation in Theorem 6.1, page 173, in \cite{RES06}.
Additional properties that we shall
use in our development are discussed in the \hyperref[SectionAppendixProperties]{Appendix}.

Define $S_{n}=X_{1}+\cdots+X_{n}+S_{0}$ for $n\geq1$. Throughout the rest
of the
paper, we shall use the notation $P_{s} ( \cdot ) $ for the
probability measure on the path-space of the process $S=(S_{n}\dvt n\geq
0)$ given
that $S_{0}=s$. Let $v_{1}^{\ast},\ldots,v_{m}^{\ast}\in\mathbb
{R}^{d}$ and
$a_{1}^{\ast},\ldots,a_{m}^{\ast}\in\mathbb{R}^{+}$. We define\vspace*{-2.4pt}
%
\begin{equation}
A^{\ast}=\bigcup_{j=1}^{m}\bigl\{y
\dvt y^{T}v_{j}^{\ast}>a_{j}^{\ast}
\bigr\}=\Bigl\{ y\dvt \max_{j=1}%
^{m}
\bigl(y^{T}v_{j}^{\ast}-a_{j}^{\ast}
\bigr)>0\Bigr\}. \label{DefAstar}%
\end{equation}
We set $T_{A^{\ast}}=\inf\{n\geq0\dvt S_{n}\in A^{\ast}\}$ and write
$bA^{\ast
}=\{y\dvt y=bx,$ $x\in A^{\ast}\}$. Note that\vspace*{-2.4pt}%
\[
bA^{\ast}=\bigl\{z\dvt r_{b}^{\ast} ( z ) >0\bigr\},
\]
where\vspace*{-2.4pt}%
\[
r_{b}^{\ast} ( z ) \triangleq\max_{j=1}^{m}
\bigl(z^{T}v_{j}^{\ast}%
-a_{j}^{\ast}b
\bigr).
\]
Finally, put
\[
u_{b}^{\ast} ( s ) =P_{s} ( T_{bA^{\ast}}<\infty
) .
\]
We are concerned with the asymptotic conditional distribution of $ (
S_{n}\dvt n\leq T_{bA^{\ast}} ) $ given that $T_{bA^{\ast}}<\infty
$ as
$b\nearrow\infty$. Throughout this paper, we impose the following two
assumptions.

\begin{assumption}\label{ass1}
$X_{n}$ has a continuous regularly varying distribution
with index $\alpha>1$ and $EX_{n}=\eta=-\mathbf{1}$, where $\mathbf
{1}= (
1,\ldots,1 ) ^{T}\in\mathbb{R}^{d}$.
\end{assumption}

\begin{assumption}\label{ass2}
For each $j$, $\eta^{T}v_{j}^{\ast}=-1$ and
$\mu (
A^{\ast} ) >0$.
\end{assumption}


\begin{remark}
Assumption \ref{ass1} indicates that $X_{n}$ has a continuous distribution. This
can be
dispensed by applying a smoothing kernel in the definition of the function
$H ( \cdot ) $ introduced later. The assumption that $\eta
=-\mathbf{1}$ is equivalent (up to a rotation and a multiplication of a
scaling factor) to $EX_{n}=\eta\neq0$, so nothing has been lost by imposing
this condition. We also assume that $\eta^{T}v_{j}^{\ast}=-1$ for
each $j$;
this, again, can always be achieved without altering the problem
structure by
multiplying the vector $v_{j}^{\ast}$ and $a_{j}^{\ast}$ by a
positive factor
as long as $\eta^{T}v_{j}^{\ast}<0$. Now, given that the random walk
has drift
$\eta$, it is not difficult to see geometrically that some conditions
must be
imposed on the $v_{i}^{\ast}$'s in order to have a meaningful large deviations
situation (i.e., $u^{*}_{b} ( 0 ) \rightarrow0$ as
$b\rightarrow
\infty$). In particular, we must have that $A^{\ast}$ does not
intersect the
ray $\{t\eta\dvt t>0\}$. Otherwise, the Law of Large Numbers might
eventually let
the process hit the set $bA^{\ast}$. However, avoiding intersection
with the
ray $\{t\eta\dvt t>0\}$ is not enough to rule out some degenerate
situations. For
instance, suppose that $A^{\ast}=\{y\dvt y^{T}v^{\ast}>1\}$ with $\eta
^{T}v^{\ast
}=0$ (i.e., the face of $A^{\ast}\ $is parallel to $\eta$); in this case
Central Limit Theorem-type fluctuations might eventually make the
random walk
hit the target set. Therefore, in order to rule out these types of degenerate
situations one requires $\eta^{T}v_{j}^{\ast}<0$.
\end{remark}

As mentioned earlier, we are interested in providing a tractable asymptotic
description of the conditional distribution of $ ( S_{n}\dvt n\leq
T_{bA^{\ast}} ) $ given that $T_{bA^{\ast}}<\infty$ as
$b\nearrow\infty
$. It is well known that, given $T_{bA^{\ast}}<\infty$, the process
$ (
S_{n}\dvt n\leq T_{bA^{\ast}} ) $ is Markovian with transition
kernel given
by%
\[
K^{\ast} ( s_{0},\mathrm{d}s_{1} ) =P ( X+s_{0}
\in \mathrm{d}s_{1} ) \frac{u_{b}^{\ast} ( s_{1} ) }{u_{b}^{\ast} (
s_{0} )
}.
\]
Note that $K^{\ast} ( \cdot ) $ is a well defined Markov transition
kernel because of the (harmonic) relationship
\[
0<u_{b}^{\ast} ( s ) =E_{s}\bigl[P_{s} (
T_{bA^{\ast
}}<\infty |X_{1} ) \bigr]=E\bigl[u_{b}^{\ast}
( s+X_{1} ) \bigr].
\]
The transition kernel $K^{\ast} ( \cdot ) $ is the Doob's
$h$-transform of the original random walk kernel (the name
$h$-transform is
given after the harmonic property of the positive function
$u_{b}^{*} (
\cdot ) $). The Markov kernel $K^{\ast} ( \cdot ) $
generates
a measure $P_{s}^{\ast} ( \cdot ) $ on the $\sigma$-field
$\mathcal{F}_{T_{bA^{\ast}}}$ generated by the $X_{k}$'s up to time
$T_{bA^{\ast}}$.

Our goal is to construct a \textit{tractable} measure $\hat{P}_{s}$ on
$\mathcal{F}_{T_{bA^{\ast}}}$ such that for each fixed $s$,%
%
\begin{equation}
\lim_{b\rightarrow\infty}\sup_{B\in\mathcal{F}_{T_{bA^{\ast
}}}}\bigl|\hat{P}%
_{s} ( B ) -P_{s}^{\ast} ( B ) \bigr|=0. \label{TV1}
\end{equation}

\begin{remark}
The measures $\hat P_{s}$ and $P^{*}_{s}$ certainly depend on the specific
rarity parameter $b$. To simplify notation, we omit the index $b$ in the
notation of $\hat P_{s}$ and $P^{*}_{s}$ when it does not cause confusion.
\end{remark}

By tractability we mean that the $\hat{P}_{s} ( \cdot ) $ is
constructed by a Markov transition kernel that can, in principle, be computed
directly from the increment distribution of the random walk. The transition
kernel associated to $\hat{P}_{s} ( \cdot ) $ will be relatively
easy to manipulate and it will be explicit in terms of the associated
increment distribution. Together with the strong mode of convergence implied
in (\ref{TV1}), we will be able to provide refinements to common
results that
are found in the literature concerning conditional limit theorems of specific
quantities of interest (cf. \cite{AsmKlup96,HultLinMikSam05}).

\subsection{Elements of the approximations and main results}

A natural strategy that one might pursue in constructing $\hat
{P}_{s} (
\cdot ) $ consists in taking advantage of the approximations
that are
available for $u_{b}^{\ast} ( s ) $.
Define%
%
\begin{eqnarray}
\label{vb} v_{b}^{\ast} ( s ) & =&\int_{0}^{\infty}P_{s}
\bigl( X+t\eta+s\in bA^{\ast} \bigr) \,\mathrm{d} t
\nonumber
\\
& =&\int_{0}^{\infty}P_{s} \Bigl( \max
_{i=1}^{m}\bigl[(s+X+t\eta )^{T}v_{i}^{\ast
}-a_{i}^{\ast}b
\bigr]>0 \Bigr) \,\mathrm{d} t
\\
& =&\int_{0}^{\infty}P_{s} \bigl(
r_{b}^{\ast}(s+X)>t \bigr) \,\mathrm{d} t=E \bigl( r_{b}^{\ast}(s+X)^{+}
\bigr) .\nonumber%
\end{eqnarray}
Note that in the third equality above we have used that $\eta
^{T}v_{i}^{\ast
}=-1$. If Assumptions \ref{ass1} and \ref{ass2} are in place, it is well known that (see
\cite{HultLin06} and \cite{HultLinMikSam05})%
%
\begin{equation}
u_{b}^{\ast} ( s ) =v_{b}^{\ast} ( s )
\bigl( 1+\mathrm{o} ( 1 ) \bigr) \label{Hult}%
\end{equation}
as $b\nearrow\infty$, uniformly over $s$ in compact sets.
This result, together with the form of the $h$-transform, suggests defining
the kernel%
\[
K_{v} ( s_{0},\mathrm{d}s_{1} ) =P ( X+s_{0}
\in \mathrm{d}s_{1} ) \frac{v_{b}^{\ast} ( s_{1} ) }{w_{b}^{\ast} (
s_{0} )
},
\]
where%
\[
w_{b}^{\ast} ( s_{0} ) =Ev_{b}^{\ast}
( s_{0}+X )
\]
is introduced to make $K_{v} ( \cdot ) $ a well defined Markov
transition kernel. It is reasonable to expect that $K_{v} ( \cdot
)
$ and the corresponding probability measure on the sample path space,
which we
denote by $P_{s}^{ ( v ) } ( \cdot ) $, will provide
good approximations to both $K^{\ast} ( \cdot ) $ and
$P_{s}^{\ast
} ( \cdot ) $. This approach is natural and it has been
successfully applied in the one dimensional setting in the context of
subexponential increment distributions in \cite{BG07AAP}. However, in the
multidimensional setting it is not entirely straightforward to evaluate and
manipulate either $v_{b}^{\ast} ( \cdot ) $ or
$w_{b}^{\ast} (
\cdot ) $. Therefore, we shall follow a somewhat different approach.

Our strategy is inspired by the way in which ruin is intuitively
expected to
occur in the context of heavy-tailed increments; namely, the underlying random
walk proceeds according to its nominal dynamics and all of a sudden a large
jump occurs which causes ruin. This intuition is made more precise by
the form
of the kernel%
%
\begin{eqnarray}\label{Dka}
\tilde{K} ( s_{0},\mathrm{d}s_{1} ) & =&p_{b} (
s_{0} ) P ( X+s_{0}\in \mathrm{d}s_{1} )
\frac{I ( X+s_{0}\in bA^{\ast} )
}{P ( X+s_{0}\in bA^{\ast} ) }+ \bigl( 1-p_{b} ( s_{0} ) \bigr) P (
X+s_{0}\in \mathrm{d}s_{1} )
\nonumber
\\[-8pt]
\\[-8pt]
\nonumber
& =&P ( X+s_{0}\in \mathrm{d}s_{1} ) \biggl\{ p_{b} (
s_{0} ) \frac{I ( X+s_{0}\in bA^{\ast} ) }{P ( X+s_{0}\in
bA^{\ast
} ) }+ \bigl( 1-p_{b} (
s_{0} ) \bigr) \biggr\} ,
\nonumber
\end{eqnarray}
where $p_{b} ( s_{0} ) $ will be suitably chosen so that
%
\begin{equation}
p_{b} ( s_{0} ) \approx P_{s_{0}} (
X_{1}+s_{0}\in bA|T_{bA^{\ast}}<\infty ) .
\label{EQpb}%
\end{equation}
In other words, $\tilde{K} ( \cdot ) $ is a mixture
involving both
the ruinous and the regular components. The mixture probability is
chosen to
capture the appropriate contribution of the ruinous component at every step.

We shall construct $\hat{P}_{s} ( \cdot ) $ by studying a
family of
transition kernels $\hat{K} ( \cdot ) $ that are very
close to~$\tilde{K} ( \cdot ) $. We will not work directly with
$\tilde
{K} ( \cdot ) $\ to avoid some uniform integrability issues that
arise in testing the Lyapunov bound to be described later in Lemma
\ref{LemLB}. The definition of $\hat{K} ( \cdot ) $
requires a
modification of the target set. This modification will be convenient because
of two reasons: first, to localize the analysis of our Lyapunov
functions only
in a suitable compact region that scales according to the parameter $b$;
second, to apply a Taylor expansion in combination with the dominated
convergence theorem. The Taylor expansion will be applied to a mollified
version of the function $r_{b}^{\ast} ( s ) $ in the verification
of the Lyapunov bound.

\subsubsection{Enlargement procedure of the target region}

First, given any $\delta\in ( 0,1 ) $ we define
$v_{j}^{\ast
} ( \delta ) =(v_{j}^{\ast}+\delta\eta/\Vert
\eta\Vert _{2}^{2})/(1-\delta)$, and observe that
$\eta
^{T}v_{j}^{\ast} ( \delta ) =-1$. We then write, given
$\beta>0$,
%
\begin{equation}
\label{A}A=A^{\ast}\cup\Biggl(\bigcup_{j=1}^{m}
\bigl\{y\dvt y^{T}v^{*}_{j}(\delta
)>a_{j}^{\ast
}\bigr\}\Biggr)\cup\Biggl(\bigcup
_{i=1}^{d}\{y\dvt y_{i}\geq\beta\}\Biggr).
\end{equation}
To simplify the notation let $e_{i}\in\mathbb{R}^{d}$ be the vector whose
$i$th component is equal to one and the rest of the components are
zero, and
express $A$ in the same form as we do for $A^{\ast}$. We write%
\[
v_{j}=\cases{ %
v_{j}^{\ast}, & \quad $1\leq j\leq m,$
\vspace*{2pt}\cr
v_{j}^{\ast} ( \delta ), & \quad $m+1\leq j\leq2m,$
\vspace*{2pt}\cr
e_{i}, &\quad  $2m+1\leq j\leq2m+d$,}
\qquad
a_{j}=\cases{a_{j}^{\ast}, & \quad $1\leq j\leq m,$
\vspace*{2pt}\cr
a_{j-m}^{\ast}, &\quad $m+1\leq j\leq2m,$
\vspace*{2pt}\cr
\beta,&\quad $2m+1\leq j\leq2m+d.$}
\]
We then have that $A=\bigcup_{j=1}^{2m+d}\{y\dvt y^{T}v_{j}>a_{j}\}$. Analogous
approximations such as (\ref{vb}) and (\ref{Hult}) are applicable. The
addition of the vectors $v_{j}$ for $j\geq m+1$ will be convenient in
order to
analyze a certain Lyapunov inequality in a compact set. Now, note that if
$u_{b} ( s ) =P ( T_{bA}<\infty ) $, then%
%
\begin{equation}
u_{b} ( s ) \bigl(1+\mathrm{o} ( 1 ) \bigr)=v_{b} ( s ) \triangleq
\int_{0}^{\infty}P ( X+t\eta+s\in bA ) \,\mathrm{d} t.
\label{Hultb}%
\end{equation}
Moreover, note that%
\[
u_{b}^{\ast} ( s ) \leq u_{b} ( s )
\]
and that, for each fixed $s$,%
%
\begin{equation}
\lim_{\delta\rightarrow0}\lim_{\beta\rightarrow\infty}\lim
_{b\rightarrow
\infty}\frac{u_{b}^{\ast} ( s ) }{u_{b} ( s ) }=1. \label{APPEND}%
\end{equation}
Our strategy consists in first obtaining results
for the event that $T_{bA}<\infty$. Then, thanks to (\ref{APPEND}),
we can
select $\beta$ arbitrarily large and $\delta$ arbitrarily small to
obtain our
stated results for the conditional distribution of the walk given
$T_{bA^{\ast}}<\infty$, which is our event of interest.

Now, define%
\[
r_{b} ( z ) \triangleq\max_{j=1}^{2m+d}
\bigl\{ \bigl(z^{T}v_{j}-a_{j}b\bigr)\bigr\},
\]
and just as we obtained for (\ref{vb}), we can conclude that%
%
\begin{equation}
v_{b} ( s ) =E \bigl( r_{b}(s+X)^{+} \bigr) .
\label{vbd}%
\end{equation}

Given $a\in ( 0,1 ) $ the enlarged region takes the form%
%
\begin{equation}
\label{Aa}A_{b,a} ( s_{0} ) =\Bigl\{s_{1}\dvt
\max_{j=1}^{2m+d}\bigl[ ( s_{1}-s_{0}
) ^{T}v_{j}-a\bigl(a_{j}b-s_{0}^{T}v_{j}
\bigr)\bigr]>0\Bigr\}.
\end{equation}
%

\subsubsection{\texorpdfstring{The family of transition kernels $\hat{K}(\cdot)$}
{The family of transition kernels K(.)}}

We now describe our proposed approximating kernel $\hat{K} (
\cdot ) $ based on the enlarged target region. Given $a\in
(0,1)$ we put%
%
\begin{equation}
\hat{K} ( s_{0},\mathrm{d}s_{1} ) =P ( X+s_{0}\in
\mathrm{d}s_{1} ) \biggl\{ \frac{p_{b} ( s_{0} ) I ( s_{1}\in
A_{b,a} (
s_{0} ) >0 ) }{P ( s_{0}+X\in A_{b,a} (
s_{0} )
>0 ) }+ \bigl( 1-p_{b} (
s_{0} ) \bigr) \biggr\} . \label{DKt}%
\end{equation}

The final ingredient in the description of our approximating kernel, and
therefore of our approximating probability measure $\hat{P}_{s} (
\cdot ) $, corresponds to the precise description of $p_{b} (
s_{0} ) $ and the specification of $a\in ( 0,1 ) $. The
scalar $a$ eventually will be chosen arbitrarily close to 1. As we indicated
before, in order to follow the intuitive description of the most likely
way in
which large deviations occur in heavy-tailed settings, we should guide the
selection of $p_{b} ( s_{0} ) $ via (\ref{EQpb}).

Using (\ref{Hultb}) and our considerations about $v_{b} (
s_{0} )
$, we have that as $b\rightarrow\infty$%
%
\begin{equation}
P_{s_{0}} ( X_{1}+s_{0}\in bA|T_{bA}<
\infty ) =\frac
{P (
r_{b} ( s_{0}+X ) >0 ) }{v_{b} ( s_{0} )
} \bigl( 1+\mathrm{o} ( 1 ) \bigr) . \label{D-1}%
\end{equation}
The underlying approximation (\ref{D-1}) deteriorates when $r_{b} (
s_{0} ) $ is not too big or, equivalently, as $s_{0}$ approaches the
target set $bA$. In this situation, it is not unlikely that the nominal
dynamics will make the process hit $bA$. Therefore, when $s_{0}$ is close
enough to $bA$, we would prefer to select $p_{b} ( s_{0} )
\approx0$. Due to these considerations and given the form of the jump set
specified in $\hat{K} ( \cdot ) $ above\ we define%
%
\begin{equation}
p_{b} ( s ) =\min \biggl( \frac{\theta P ( s_{0}+X\in
A_{b,a} ( s_{0} )  ) }{v_{b} ( s_{0} )
},1 \biggr) I \bigl(
r_{b} ( s ) \leq-\delta_{2}b \bigr) \label {D-p}%
\end{equation}
for $\delta_{2}>0$ chosen small enough and $\theta,a\in (
0,1 ) $
chosen close enough to 1. The selection of all these constants will be
done in
our development.

\subsubsection{The statement of our main results}

Before we state our main result, we need to describe the probability
measure in
path space that we will use to approximate
\[
P_{s}^{\ast}(S\in\cdot)\triangleq P(S\in\cdot|T_{bA^{\ast
}}<
\infty,S_{0}=s),
\]
in total variation. Given $\gamma>0$ define%
\[
\Gamma=\bigl\{y\dvt y^{T}\eta\geq\gamma\bigr\}
\]
and set $T_{bA}=\inf\{n\geq0\dvt S_{n}\in bA\}$, and $T_{b\Gamma}=\inf
\{n\geq0\dvt S_{n}\in b\Gamma\}$. We now define the change of measure
that we
shall use to approximate $P_{s}^{\ast}(\cdot)$ in total variation.
\vspace*{-2pt}

\begin{definition}\label{def1}
Let $\hat{P}_{s} ( \cdot ) $ be
defined as the
measure generated by transitions according to $\hat{K} ( \cdot
) $
up to time $T_{bA}\wedge T_{b\Gamma}$, with mixture probability as
defined in
\eqref{D-p}, and transitions according to $K ( \cdot ) $
for the
increments $T_{bA}\wedge T_{b\Gamma}+1$ up to infinity.\vspace*{-2pt}
\end{definition}

We now state our main result.\vspace*{-2pt}

\begin{theorem}
\label{ThmMain1}For every $\varepsilon>0$ there exists $\theta
,a,\delta
,\delta_{2}\in ( 0,1 ) $ ($\theta,a$ sufficiently close
to 1 and
$\delta,\delta_{2}$ sufficiently close to zero), and $\beta,\gamma,b_{0}>0$
sufficiently large so that if $b\geq b_{0}$
\[
\sup_{B\in\mathcal{F}}\bigl|\hat{P}_{0} ( B ) -P_{0}^{\ast
}
( B ) \bigr|\leq\varepsilon,
\]
where $\mathcal{F}=\sigma(\bigcup_{n=0}^{\infty}\sigma\{S_{k}\dvt 0\leq
k\leq n\})$.\vspace*{-2pt}
\end{theorem}

In order to illustrate an application of the previous result, we have
the next
theorem which follows without much additional effort, as a corollary to
Theorem \ref{ThmMain1}. The statement of the theorem, however,
requires some
definitions that we now present.

Because of regular variation, we can define for any $a_{1}^{\ast}%
,\ldots,a_{m}^{\ast}>0$
%
\begin{equation}
\lim_{b\rightarrow\infty}\frac{P ( \max_{j=1}^{m}[X^{T}v_{j}^{\ast}%
-a_{j}^{\ast}b]>0 ) }{P ( \Vert  X\Vert _{2}>b ) }=\kappa^{\ast}
\bigl(a_{1}^{\ast
},\ldots,a_{m}^{\ast}\bigr),
\label{kappa}%
\end{equation}
where for any $t\geq0$%
\[
\kappa\bigl(a_{1}^{\ast}+t,\ldots,a_{m}^{\ast}+t
\bigr)\triangleq\mu \Bigl( \Bigl\{ y\dvt \max_{j=1}^{m}
\bigl( y^{T}v_{j}^{\ast}-a_{j}^{\ast}
\bigr) >t\Bigr\} \Bigr) .
\]
Using this representation, we obtain that
%
\begin{eqnarray}
\label{AvS} v_{b}^{\ast}(s) & =&\int_{0}^{\infty}P
\Bigl( \max_{j=1}^{m}\bigl[(s+X)^{T}%
v_{j}^{\ast}-a_{j}^{\ast}b\bigr]>t \Bigr) \,\mathrm{d} t
\nonumber
\\
& =& \bigl( 1+\mathrm{o}  ( 1 ) \bigr) P \bigl( \Vert X\Vert
_{2}>b \bigr)\nonumber\\
 &&{}\times \int_{0}^{\infty}\kappa
\bigl( a_{1}^{\ast}-b^{-1}s^{T}v_{1}^{\ast}+b^{-1}t,\ldots,a_{m}^{\ast
}-b^{-1}s^{T}%
v_{m}^{\ast}+b^{-1}t \bigr) \,\mathrm{d} t
\\
& =& \bigl( 1+\mathrm{o}  ( 1 ) \bigr) bP \bigl( \Vert X\Vert
_{2}>b \bigr) \int_{0}^{\infty}\kappa
\bigl( a_{1}^{\ast}-b^{-1}s^{T}v_{1}^{\ast}+t,\ldots,a_{m}^{\ast
}-b^{-1}s^{T}v_{m}%
^{\ast}+t \bigr) \,\mathrm{d} t\nonumber
\\
& =& \bigl( 1+\mathrm{o}  ( 1 ) \bigr) bP \bigl( \Vert X\Vert
_{2}>b \bigr) \int_{0}^{\infty}\kappa
\bigl( a_{1}^{\ast}+t,\ldots,a_{m}^{\ast}+t \bigr)
\,\mathrm{d} t\nonumber %
\end{eqnarray}
as $b\rightarrow\infty$ uniformly over $s$ in a compact set.
Actually, an
extension to approximation (\ref{AvS}) to the case in which $s=\mathrm{O}  (
b ) $ is given in the \hyperref[SectionAppendixProperties]{Appendix}; see Lemma \ref{LemRatv}. To further
simplify the notation, we write%
%
\begin{equation}
\kappa_{\mathbf{a}^{\ast}}(t)=\kappa\bigl(a_{1}^{\ast
}+t,\ldots,a_{m}^{\ast}+t
\bigr), \label{kappaa}%
\end{equation}
where $\mathbf{a}^{\ast}=(a_{1}^{\ast},\ldots,a_{m}^{\ast})$.

\begin{theorem}
\label{ThmMain2}For each $z>0$ let $Y^{\ast} ( z ) $ be a random
variable with distribution given by
\[
P \bigl( Y^{\ast} ( z ) \in B \bigr) =\frac{\mu (
B\cap\{y\dvt \max_{j=1}^{m}[y^{T}v_{j}^{\ast}-a_{j}^{\ast}]\geq z\}
) }%
{\mu
( \{y\dvt \max_{j=1}^{m}
[y^{T}v_{j}^{\ast}-a_{j}^{\ast}
]\geq z\} ) }.
\]
In addition, let $Z^{\ast}$ be a positive random variable following
distribution
\[
P \bigl( Z^{\ast}>t \bigr) =\exp \biggl\{ -\int_{0}^{t}
\frac{\kappa
_{\mathbf{a}^{\ast}}(s)}{\int_{s}^{\infty}\kappa_{\mathbf
{a}^{\ast}}%
(u)\,\mathrm{d} u}\,\mathrm{d} s \biggr\}
\]
for $t\geq0$ where $\kappa_{\mathbf{a}^{\ast}} ( \cdot
) $ is as
defined in (\ref{kappaa}). Then if $S_{0}=0$ and $\alpha>2$, we have
that%
\[
\biggl( \frac{T_{bA^{\ast}}}{b},\frac{S_{uT_{bA^{\ast}}}-u
T_{bA^{\ast}}\eta
}{\sqrt{T_{bA^{*}}}},\frac{X_{T_{bA^{\ast}}}}{b} \biggr)
\Rightarrow \bigl( Z^{\ast},CB \bigl( uZ^{\ast} \bigr)
,Y^{\ast} \bigl( Z^{\ast
} \bigr) \bigr)
\]
in $\mathbb{R}\times D[0,1)\times\mathbb{R}^{d}$, where $CC^{T}=\operatorname{Var}(X)$,
$B ( \cdot ) $ is a $d$-dimensional Brownian motion with identity
covariance matrix, $B(\cdot)$ is independent of $Z^{\ast}$ and
$Y^{\ast
} ( Z^{\ast} ) $.
\end{theorem}

\begin{remark}\label{rem3}
The random variable $Z^{\ast}$ (multiplied by a factor of $b$)
corresponds to
the asymptotic time to ruin. In the one dimensional setting, $Z^{\ast}$
follows a Pareto distribution with index $\alpha-1$. The reason is
that in the
one dimensional case%
\[
\kappa_{\mathbf{a}^{\ast}}(s)=\frac{\alpha-1}{s}\int_{s}^{\infty}
\kappa_{\mathbf{a}^{\ast}}(u)\,\mathrm{d} u.
\]
This no longer can be ensured in the multidimensional case. Nevertheless,
$Z^{\ast}$ is still regularly varying with index $\alpha-1$.
\end{remark}

\begin{remark}\label{rem4}
If $\alpha\in(1,2]$, then our analysis allows to conclude that%
\[
\biggl( \frac{T_{bA^{\ast}}}{b},\frac{S_{uT_{bA^{\ast
}}}}{T_{bA^{\ast}}}%
,\frac{X_{T_{bA^{\ast}}}}{b}
\biggr) \Rightarrow \bigl( Z^{\ast
},u\eta ,Y^{\ast} \bigl(
Z^{\ast} \bigr) \bigr)
\]
in $\mathbb{R}\times D[0,1)\times\mathbb{R}^{d}$ as $b\rightarrow
\infty$.\vadjust{\goodbreak}
\end{remark}

\section{Total variation approximations and Lyapunov inequalities}\label{SecLBTV}

In this section, we provide the proof of Theorem \ref{ThmMain1}.
First, it is
useful to summarize some of the notation that has been introduced so far.

\begin{enumerate}
\item The set $A^{\ast}$ be the target set, $v_{b}^{\ast}(s)$ be the
approximation of $P_{s}(T_{bA^{\ast}}<\infty)$, and $P_{s}^{\ast}%
(\cdot)=P(\cdot|T_{bA^{\ast}}<\infty)$ be the corresponding
conditional distribution.

\item The set $A$ is an enlargement of $A^{\ast}$ and depends on
$\delta$ and
$\beta$; $v_{b}(s)$ be the approximation of $P_{s}(T_{bA}<\infty)$.

\item The set $\Gamma=\{y\dvt y^{T}\eta\geq\gamma\}$ will be used to
define an
auxiliary conditional distribution below.

\item Under the distribution $\hat{P}_{s}(\cdot)$ in path space increments
follow the transition kernel $\hat{K}$ up to time $T_{bA}\wedge
T_{b\Gamma}$.
\end{enumerate}

Now, we shall outline the program that will allow us to proof Theorem
\ref{ThmMain1}. The program contains three parts. The first part
consists in
introducing an auxiliary conditional distribution involving a finite horizon.
To this end, we define%
%
\begin{equation}
P_{s}^{\&} ( \cdot ) \triangleq P(S\in\cdot|T_{bA}
\leq T_{b\Gamma
},S_{0}=s). \label{DefPPP}%
\end{equation}
Eventually, as we shall explain, we will select $\delta>0$ is small enough,
$\beta$ and $\gamma$ are large enough. The second part consists in showing
that $P_{s}^{\&}$ and $P_{s}^{\ast}$ are close in total variation;
this will
be done in Lemma \ref{LemFiniteHorizon}. Finally, in the third part we show
that $P_{s}^{\&}$ can be approximated by $\hat{P}_{s}$ in total variation;
this will be done in Proposition \ref{PropFinieH}. Theorem 1 then follows
directly by combining Lemma \ref{LemFiniteHorizon} and Proposition
\ref{PropFinieH}.

In order to carry out the third part of our program, namely, approximating
$P_{s}^{\&}$ by $\hat{P}_{s}$ in total variation. A natural approach,
which we
shall follow, is to argue that $\mathrm{d}P_{s}^{\&}/\mathrm{d}\hat{P}_{s}$ is close to unity.
We define
\[
\beta ( s ) \triangleq\hat{E}_{s}\bigl((\mathrm{d}P_{s}/\mathrm{d}\hat
{P}_{s})^{2}I ( T_{bA}\leq T_{b\Gamma} )
\bigr)=P_{s} ( T_{bA}\leq T_{b\Gamma
} ) ^{2}
\hat{E}_{s}\bigl(\bigl(\mathrm{d}P_{s}^{\&}/\mathrm{d}
\hat{P}_{s}\bigr)^{2}\bigr),
\]
note that $\beta ( s ) \geq P_{s} ( T_{bA}\leq
T_{b\Gamma
} ) ^{2}$ (by Jensen's inequality). As we shall verify in Lemma
\ref{LemTV1}, if we are able to show that $\beta ( 0 )
\leq
P_{0} ( T_{bA}\leq T_{b\Gamma} ) ^{2}(1+\varepsilon)$
(thus showing
that $\mathrm{d}P_{s}^{\&}/\mathrm{d}\hat{P}_{s}$ is close to unity) then we will be
able to
claim that $P_{0}^{\&}$ and $\hat{P}_{0}$ are close in total variation.

Obtaining a useful bound for $\beta ( s ) $ is the most demanding
part of the whole program. The strategy relies on the so-called Lyapunov
inequalities. The idea is to find a function $g ( \cdot )
$, which
is called a Lyapunov function, satisfying certain criteria specified in Lemma
\ref{LemLB} in order to ensure that $g ( s ) \geq\beta
(
s ) $. Now, constructing Lyapunov functions is not easy, however, we
eventually wish to enforce an upper bound corresponding to the behavior
of $P_{s} ( T_{bA}\leq T_{b\Gamma} ) ^{2}$, so it makes
sense to
use $v_{b} ( \cdot ) ^{2}$ (recall equation (\ref{Hult}) and
(\ref{vbd})) as starting template for $g ( \cdot ) $. In the
process of verifying that criteria in Lemma \ref{LemLB} it is useful
to ensure
some smoothness properties of a candidate Lyapunov function. So, a
mollification procedure is performed to the template suggested by
$v_{b} ( \cdot ) ^{2}$. The verification of the criteria in Lemma
\ref{LemLB} is pursued in Section \ref{SubsectionVerification}.

Now we start executing the first part of the previous program. Before
we start
it is useful to remark that in our definition of $P_{s}^{\&} (
\cdot ) $ (see (\ref{DefPPP})),\ we write $T_{bA}\leq
T_{b\Gamma}$ rather
than $T_{bA}<T_{b\Gamma}$. This distinction is important in the proof
of the
next result because, due to the geometry of the sets $A$ and $A^{\ast
}$, on
the set $T_{bA}>T_{b\Gamma}$ we can guarantee that $S_{T_{b\Gamma}}$ is
sufficiently far away from the set $bA^{\ast}$.

\begin{lemma}
\label{LemFiniteHorizon}For each $\varepsilon>0$ we can find $\delta
,\beta,\gamma>0$ such that%
\[
\mathop{\overline{\lim}}_{b\rightarrow\infty}\biggl\llvert \frac{P_{0} (
T_{bA^{\ast
}}<\infty ) }{P_{0} ( T_{bA}\leq T_{b\Gamma} )
}-1\biggr\rrvert
\leq\varepsilon;
\]
moreover, for $b$ sufficiently large
\[
\sup_{B\in\mathcal{F}}\bigl\llvert P_{0}^{\ast} ( B )
-P_{0}%
^{\&} ( B ) \bigr\rrvert \leq\varepsilon.
\]

\end{lemma}

\begin{pf}
We prove the first part of the lemma by establishing an upper and lower bound
of $P_{0}(T_{bA^{*}}<\infty),$ respectively.

\textit{Upper bound}.
Observe that%
\begin{eqnarray*}
P_{0} ( T_{bA^{\ast}}<\infty ) & =&P_{0} (
T_{bA^{\ast
}%
}<\infty,T_{bA}\leq T_{b\Gamma} ) +P_{0}
( T_{bA^{\ast}}%
<\infty,T_{bA}>T_{b\Gamma} )
\\
&    \leq& P_{0} ( T_{bA}\leq T_{b\Gamma} ) +\sup
_{\{s\dvt s\in
b\Gamma,s\notin bA\}}P_{s} ( T_{bA^{\ast}}<\infty ) .
\end{eqnarray*}
Note that if $s\in b\Gamma$ and $s\notin bA$ then $s^{T}\eta\geq
\gamma b$,
$s^{T}v_{i}=s^{T}v_{i}^{\ast}+\delta s^{T}\eta/\Vert
\eta\Vert _{2}^{2}\leq a_{i}b ( 1-\delta
) $ for
$i\in\{1,\ldots,m\}$. Therefore,
\[
s^{T}v_{i}^{\ast}\leq a_{i}b ( 1-\delta )
-\delta s^{T}%
\eta/\Vert \eta\Vert _{2}^{2}\leq a_{i}b-\delta\bigl(
\gamma+a_{i} \Vert \eta\Vert ^{2}_{2}\bigr) b/
\Vert\eta \Vert _{2}^{2}.
\]
Thus, $\sup_{s\in b\Gamma\setminus bA} s^{T} v^{*}_{i} - a_{i} b
\leq
-\delta(\gamma+a_{i}\Vert \eta\Vert ^{2}_{2})b/\Vert \eta\Vert ^{2}_{2}$. Given
$\varepsilon
>0$, after choosing $\delta>0$, we can select $\gamma$ large enough
so that
\[
\sup_{\{s\dvt s\in b\Gamma,s\notin bA\}}P_{s}( T_{bA^{\ast}}<\infty) \leq
\varepsilon P_{0} ( T_{bA^{\ast}}<\infty ). 
\]
%
Therefore, we have that given $\varepsilon>0$, we can select $\delta>0$
sufficiently small and $\gamma$ large enough so that for all $b$ large
enough%
\[
P_{0} ( T_{bA^{\ast}}<\infty ) \leq P_{0} (
T_{bA}\leq T_{b\Gamma} ) +\varepsilon P_{0} (
T_{bA^{\ast}}<\infty ) ,
\]
which yields an upper bound of the form
\[
P_{0} ( T_{bA^{\ast}}<\infty ) ( 1-\varepsilon ) \leq
P_{0} ( T_{bA}\leq T_{b\Gamma} ) .
\]

\textit{Lower bound}.
Notice that $P(T_{b\Gamma}<\infty)=1$. Then, we have that%
%
\begin{equation}
P_{0} ( T_{bA}\leq T_{b\Gamma} ) \leq P_{0}
( T_{bA}%
<\infty ) =v_{b} ( 0 ) \bigl(1+\mathrm{o}  ( 1 )
\bigr) \label{BndTg}%
\end{equation}
as $b\rightarrow\infty$. In addition, by the asymptotic approximation
of the
first passage time probability, Lemmas \ref{LemPropmu} and \ref
{LemRatv} in
the \hyperref[SectionAppendixProperties]{Appendix}, given $\varepsilon>0$ we can select $\delta,\beta>0$
such that
%
\begin{equation}
1\leq\mathop{\overline{\lim}}_{b\rightarrow\infty}\frac{v_{b} (
0 ) }%
{v_{b}^{\ast} ( 0 ) }=\mathop{\overline{\lim}}_{b\rightarrow
\infty}%
\frac{P_{0} ( T_{bA}<\infty ) }{P_{0} ( T_{bA^{*}}
<\infty ) }\leq1+\varepsilon. \label{Ratvs}%
\end{equation}
Thus, we obtain that as $b\rightarrow\infty$
\[
P_{0} ( T_{bA}\leq T_{b\Gamma} ) \leq\bigl(1+
\varepsilon+ \mathrm{o} (1)\bigr)P_{0}%
( T_{bA^{*}} <\infty ) .
\]
We conclude the lower bound and the first part of the lemma.

The second part of the lemma follows as an easy consequence of the
first part.
\end{pf}

Throughout our development, we then concentrate on approximating in total
variation $P_{s}^{\&} ( \cdot ) $. We will first prove the
following result.

\begin{proposition}
\label{PropFinieH}For all $\varepsilon,\delta,\beta>0$ there exists
$\theta,a$ sufficiently close to $1$ from below, $\delta_{2} $ sufficiently
small, and $\gamma$, $b$ sufficiently large such that
\[
\sup_{B\in\mathcal{F}}\bigl|\hat{P}_{0} ( B ) -P_{0}^{\&
}
( B ) \bigr|\leq\varepsilon,
\]
where $\mathcal{F}=\sigma(\bigcup_{n=0}^{\infty}\sigma\{S_{k}\dvt 0\leq
k\leq n\})$.
\end{proposition}

As noted earlier Proposition \ref{PropFinieH} combined with Lemma
\ref{LemFiniteHorizon} yields the proof of Theorem \ref{ThmMain1}. To provide
the proof of Proposition \ref{PropFinieH}, we will take advantage of the
following simple yet powerful observation (see also \cite{LECBLATUFGLY}).

\begin{lemma}
\label{LemTV1}Let $Q_{0}$ and $Q_{1}$ be probability measures defined
on the
same $\sigma$-field $\mathcal{G}$ and such that $\mathrm{d}Q_{1}=M^{-1}\,\mathrm{d} Q_{0}$
for a
positive r.v. $M>0$. Suppose that for some $\varepsilon>0$,
$E^{Q_{1}} (
M^{2} ) =E^{Q_{0}}M\leq1+\varepsilon$. Then
\[
\sup_{B\in\mathcal{G}}\bigl\llvert Q_{1} ( B )
-Q_{0} ( B ) \bigr\rrvert \leq\varepsilon^{1/2}.
\]

\end{lemma}

\begin{pf}
Note that%
\begin{eqnarray*}
\bigl\llvert Q_{1} ( B ) -Q_{0} ( B ) \bigr\rrvert & =&
\bigl\llvert E^{Q_{1}} ( 1-M;B ) \bigr\rrvert
\\
& \leq& E^{Q_{1}} \bigl( \llvert M-1\rrvert \bigr) \leq E^{Q_{1}%
}
\bigl[ ( M-1 ) ^{2}\bigr]^{1/2}= \bigl( E^{Q_{1}}M^{2}-1
\bigr) ^{1/2}%
\leq\varepsilon^{1/2}.
\end{eqnarray*}
\upqed\end{pf}

Lemma \ref{LemTV1} will be used to prove Proposition \ref
{PropFinieH}. In
particular, we will first apply Lemma \ref{LemTV1} by letting
$\mathcal{G}%
=\mathcal{F}_{T_{bA}\wedge T_{b\Gamma}}$, $Q_{1}=\hat{P}_{s}$, and
$Q_{0}=P_{s}^{\&} ( \cdot ) $. The program proceeds as follows.
With $\hat{K}$ defined as in \eqref{DKt}, we let
%
\begin{eqnarray}
\label{r} \hat{k} ( s_{0},s_{1} ) & =&\frac{P ( X+s_{0}\in
\mathrm{d}s_{1} ) }{\hat{K} ( s_{0},\mathrm{d}s_{1} ) }
\nonumber
\\
& =& \biggl\{ \frac{p_{b} ( s_{0} ) I ( s_{1}\in
A_{b,a} (
s_{0} )  ) }{P ( s_{0}+X\in A_{b,a} ( s_{0} )
) }+ \bigl( 1-p_{b} ( s_{0} )
\bigr) \biggr\} ^{-1}
\nonumber
\\[-8pt]
\\[-8pt]
\nonumber
& =&\frac{P ( s_{0}+X\in A_{b,a} ( s_{0} )  )
I (
s_{1}\in A_{b,a} ( s_{0} )  ) }{p_{b} (
s_{0} )
+ ( 1-p_{b} ( s_{0} )  ) P ( s_{0}+X\in
A_{b,a} ( s_{0} )  ) }
\\
&&{} +\frac{1}{ ( 1-p_{b} ( s_{0} )  ) }I \bigl( s_{1}\notin A_{b,a} (
s_{0} ) \bigr) .\nonumber %
\end{eqnarray}
Observe that on $\mathcal{G}$ we have%
\[
\frac{\mathrm{d}P_{s}^{\&}}{\mathrm{d}\hat{P}_{s}}=\frac{I ( T_{bA}\leq
T_{b\Gamma} )
}{P_{s} ( T_{bA}\leq T_{b\Gamma} ) }\times\frac
{\mathrm{d}P_{s}}{\mathrm{d}\hat
{P}_{s}}=\frac{I ( T_{bA}\leq T_{b\Gamma} ) }{P_{s} (
T_{bA}\leq T_{b\Gamma} ) }
\prod_{j=0}^{T_{bA}-1}\hat{k} (
S_{j},S_{j+1} ) .
\]
Therefore, according to Lemma \ref{LemTV1}, it suffices to show that%
\begin{eqnarray*}
&& \hat{E}_{s} \biggl( \biggl( \frac{\mathrm{d}P_{s}^{\&}}{\mathrm{d}\hat{P}_{s}} \biggr)
^{2} \biggr)
\\
&&\quad =\frac{1}{P_{s} ( T_{bA}\leq T_{b\Gamma} ) ^{2}}\hat{E}%
_{s} \biggl( \biggl(
\frac{\mathrm{d}P_{s}}{\mathrm{d}\hat{P}_{s}} \biggr) ^{2}I ( T_{bA}\leq
T_{b\Gamma} ) \biggr)
\\
&&\quad =\frac{1}{P_{s} ( T_{bA}\leq T_{b\Gamma} )
^{2}}E_{s} \Biggl( \prod
_{j=0}^{T_{bA}-1}\hat{k} ( S_{j},S_{j+1}
) I ( T_{bA}\leq T_{b\Gamma} ) \Biggr) \leq1+\varepsilon
\end{eqnarray*}
for all $b$ sufficiently large. Consequently, we must be able to
provide a
good upper bound for the function
\[
\hat{\beta}_{b} ( s ) \triangleq E_{s} \Biggl( \prod
_{j=0}%
^{T_{bA}-1}\hat{k} (
S_{j},S_{j+1} ) I ( T_{bA}\leq T_{b\Gamma} )
\Biggr) .
\]
In order to find an upper bound for $\hat{\beta}_{b} ( s )
$ we
will construct an appropriate Lyapunov inequality based on the following
lemma, which follows as in Theorem 2 part (iii) of \cite{BG07AAP}; see also
Theorem~2.6 in \cite{konmey05a}.

\begin{lemma}
\label{LemLB} Suppose that $C$ and $B$ are given sets. Let $\tau_{B} =
\inf\{n\dvt S_{n} \in B\}$ and $\tau_{C} = \inf\{n\dvt S_{n} \in C\}$ be
the first
passage times. Assume that $g ( \cdot ) $ is a non-negative
function satisfying%
%
\begin{equation}
g ( s ) \geq E_{s} \bigl( \hat{k} ( s,S_{1} ) g (
S_{1} ) \bigr) \label{LI-1}%
\end{equation}
for $s\notin C\cup B$. Then,
\[
g ( s ) \geq E_{s} \Biggl( g ( S_{\tau_{C}} ) \prod
_{j=0}^{\tau_{C}-1}\hat{k} ( S_{j},S_{j+1}
) I ( \tau_{C}\leq\tau_{B},\tau_{C}<\infty )
\Biggr) .
\]
Furthermore, let $h ( \cdot ) $ be any non-negative function and
consider the expectation%
\[
\hat{\beta}_{b}^{h} ( s ) =E_{s} \Biggl( h (
S_{\tau
_{C}%
} ) \prod_{j=0}^{\tau_{C}-1}\hat{k} (
S_{j},S_{j+1} ) I ( \tau_{C}\leq
\tau_{B},\tau_{C}<\infty ) \Biggr) .
\]
If in addition to (\ref{LI-1}) we have that $g ( s ) \geq
h (
s ) $ for $s\in C$, then we conclude that%
\[
g ( s ) \geq\hat{\beta}_{b}^{h} ( s ) .
\]

\end{lemma}


We will construct our Lyapunov function in order to show that for any given
$\varepsilon>0$
%
\begin{equation}
\mathop{\overline{\lim}}_{b\rightarrow\infty}\frac{\hat{\beta
}_{b}^{h} ( 0 )
}{P_{0} ( T_{bA}\leq T_{b\Gamma} ) ^{2}}\leq1+\varepsilon,
\label{lim1}%
\end{equation}
with $h\approx1$. Note that given any $\varepsilon>0$ we can select
$\gamma>0$
sufficiently large so that for all $b$ large enough%
\[
P_{0} ( T_{bA}\leq T_{b\Gamma} ) \leq P_{0}
( T_{bA}%
<\infty ) \leq ( 1+\varepsilon ) P_{0} (
T_{bA}\leq T_{b\Gamma} ) .
\]
Thus, using a completely analogous line of thought leading to (\ref
{vb}) we
would like to construct a Lyapunov function
%
\begin{equation}
g_{b} ( s ) \approx v_{b}^{2} ( s ) \triangleq
E \bigl( r_{b} ( s+X ) ^{+} \bigr) ^{2}.
\label{IntuiDefg}%
\end{equation}
If such a selection of a Lyapunov function is applicable then, given
$\varepsilon>0$, if $\gamma$ is chosen sufficiently large, we would
be able to
conclude the bound (\ref{lim1}).

\subsection{Mollification of $r_{b}(s)$ and proposal of $g_{b}(s)$}\label{SubsectionMollification}

In order to verify the Lyapunov inequality from Lemma \ref{LemLB} it
will be
useful to perform a Taylor expansion of the function $g_{b} (
\cdot ) $. Since the function $r_{b} ( s+X ) ^{+}$ is not
smooth in $s$, we will first perform a mollification procedure. Given
$c_{0}>0$ define
\[
\varrho_{b} ( s ) =c_{0}\log\Biggl(\sum
_{j=1}^{2m+d}\exp\bigl(\bigl[s^{T} 
v_{j}-a_{j}b\bigr]/c_{0}\bigr)\Biggr)
\]
and note that%
%
\begin{equation}
r_{b} ( s ) \leq\varrho_{b} ( s ) \leq r_{b} (
s ) +c_{0}\log ( 2m+d ) . \label{SInrho}
\end{equation}

Then, for $\delta_{0}>0$ let
\[
d ( x ) =\cases{ %
0, &\quad  $x\leq-
\delta_{0},$
\vspace*{2pt}\cr
( x+\delta_{0} ) ^{2}/(4\delta_{0}), &\quad $\llvert x
\rrvert \leq\delta_{0},$
\vspace*{2pt}\cr
x, & \quad $x\geq\delta_{0}.$}
\]
Further note that $d ( \cdot ) $ is continuously differentiable
with derivative $d^{\prime} ( \cdot ) $ given by the function
\[
d^{\prime}(x)=\frac{x+\delta_{0}}{2\delta_{0}}I \bigl( \llvert x\rrvert \leq
\delta_{0} \bigr) +I ( x\geq\delta_{0} ) \leq I ( x\geq-
\delta_{0} )
\]
and that
%
\begin{equation}
x^{+}\leq d ( x ) \leq ( x+\delta_{0} ) ^{+}.
\label{Ind}%
\end{equation}
Note that the functions $\varrho_{b} ( \cdot ) $ and
$d (
\cdot ) $ depend on $c_{0}$ and $\delta_{0}$, respectively; and
that we
have chosen to drop this dependence in our notation. The selection of
$\delta_{0}$ is quite flexible given that we will eventually send
$b\rightarrow\infty$. We choose
%
\begin{equation}
c_{0}\triangleq c_{0}( b) =\max\bigl(b^{{(3-\alpha)}/{2}},
\tilde{c} _{0}\bigr) \label{c0}%
\end{equation}
for some constant $\tilde{c}_{0}>0$ chosen sufficiently large. Note
that we
have
\[
b^{2-\alpha}=\mathrm{o} \bigl(c_{0} ( b ) \bigr),\qquad  c_{0} ( b )
=\mathrm{o} (b).
\]
Given $c_{0}$ and $\delta_{0}$ first we define
%
\begin{equation}
H_{b} ( s ) =E\bigl[d \bigl( \varrho_{b} ( s+X ) \bigr)
\bigr]. \label{DefHb}%
\end{equation}
We will see that the asymptotics of $H_{b} ( s ) $ as $b$
goes to
infinity are independent of $c_{0}$ and $\delta_{0}$. Indeed, observe, using
inequality (\ref{Ind}) that%
%
\begin{eqnarray}
\label{InR} E\bigl[\varrho_{b}^{+} ( s+X ) \bigr] & =&\int
_{0}^{\infty}P \bigl( \varrho_{b} ( s+X ) >t
\bigr) \,\mathrm{d} t\leq H_{b} ( s )
\nonumber
\\[-8pt]
\\[-8pt]
\nonumber
& \leq& E\bigl[ \bigl( \varrho_{b} ( s+X ) +\delta_{0}
\bigr) ^{+}\bigr]=\int_{0}^{\infty}P \bigl(
\varrho_{b} ( s+X ) >t-\delta _{0} \bigr) \,\mathrm{d} t. %
\end{eqnarray}
Therefore, using (\ref{SInrho}), the inequalities (\ref{InR}) and basic
properties of regularly varying functions (i.e., that regularly varying
functions possess long tails), we obtain that for any given $\delta_{0}>0$
%
\begin{eqnarray}
H_{b} ( s ) & = &\bigl( 1+\mathrm{o}  ( 1 ) \bigr) \int_{0}^{\infty}P
\bigl( \varrho_{b} ( s+X ) >t \bigr) \,\mathrm{d} t\label{UB}
\nonumber\\
& =& \bigl( 1+\mathrm{o}  ( 1 ) \bigr) \int_{0}^{\infty}P \bigl(
\varrho_{b} ( s+X ) >t-\delta_{0} \bigr) \,\mathrm{d} t
\\
& = &\bigl( 1+\mathrm{o}  ( 1 ) \bigr) \int_{0}^{\infty}P \Bigl(
\max_{i=1}^{2m+d}\bigl\{(s+X)^{T}v_{i}-a_{i}b
\bigr\}>t+\mathrm{o} (b) \Bigr) \,\mathrm{d} t
\nonumber\\
& =& \bigl( 1+\mathrm{o}  ( 1 ) \bigr) v_{b} ( s ) =\bigl(1+\mathrm{o} (1)
\bigr)P_{s} ( T_{bA}<\infty )
\nonumber
\end{eqnarray}
as $b\rightarrow\infty$ uniformly over $s$ in any compact set.

Finally,
$H_{b} ( \cdot ) $ is \textit{twice} continuously differentiable.
This is a desirable property for the verification of our Lyapunov inequality.
Therefore, \textit{we define the candidate of Lyapunov function, }%
$g_{b} ( \cdot ) $ via
\[
g_{b} ( s ) =\min\bigl(c_{1}H_{b} ( s )
^{2},1\bigr),
\]
where $c_{1}\in ( 1,\infty ) $ will be chosen arbitrarily
close to
$1$ if $b$ is sufficiently large. The intuition behind the previous selection
of $g_{b} ( s ) $ has been explained in the argument leading to
(\ref{IntuiDefg}).

\subsection{Verification of the Lyapunov inequality}\label{SubsectionVerification}

We first establish the Lyapunov inequality (\ref{LI-1}) on the region
$s\notin
b\Gamma$ and $r_{b} ( s ) \leq-\delta_{2}b$ for $\delta_{2}>0$
suitably small and $b$ sufficiently large. If $r_{b} ( s )
\leq-\delta_{2}b$ for $b$ large enough so that $g_{b} ( s
) <1$
then, using expression (\ref{r}), we note that inequality (\ref
{LI-1}) is
equivalent to%
%
\begin{equation}
J_{1}+J_{2}\leq1, \label{LI-1b}%
\end{equation}
where%
\[
J_{1}  =E \biggl( \frac{g_{b} ( s+X ) }{g_{b}(s)};s+X\in A_{b,a} ( s )
\biggr) \times\frac{P ( s+X\in
A_{b,a} (
s )  ) }{p_{b} ( s ) + ( 1-p_{b} (
s )
) P(s+X\in A_{b,a} ( s ) )}
\]
and%
\[
\bigl( 1-p_{b} ( s ) \bigr) J_{2}=E \biggl(
\frac{g_{b}%
(s+X)}{g_{b}(s)};s+X\notin A_{b,a} ( s ) \biggr) .
\]

In order to verify (\ref{LI-1b}) on the region $s\notin b\Gamma$,
$r_{b} ( s ) \leq-\delta_{2}b$ we will apply a Taylor
expansion to
the function $g_{b} ( \cdot ) $. This Taylor expansion will be
particularly useful for the analysis of $J_{2}$. The following result, which
summarizes useful properties of the derivatives of $\varrho_{b} (
\cdot ) $ and $d ( \cdot ) $, will be useful.

\begin{lemma}
\label{PropertiesRho}Let%
\[
w_{j} ( s ) =\frac{\exp (  (
s^{T}v_{j}-a_{j}b )
/c_{0} ) }{\sum_{i=1}^{2m+d}\exp (  ( s^{T}v_{i}%
-a_{i}b ) /c_{0} ) }.
\]
Then,
\begin{longlist}[(ii)]
\item[(i)] $\triangledown\varrho_{b} ( s ) =\sum_{j=1}^{2m+d}%
v_{j}w_{j} ( s ) $,

\item[(ii)] $\Delta\varrho_{b} ( s ) =\sum_{j=1}^{2m+d}%
w_{j} ( s )  ( 1-w_{j} ( s )  ) v_{j}%
v_{j}^{T}/c_{0}$.
\end{longlist}
\end{lemma}

\begin{pf}
Item (i) follows from basic calculus. Part (ii) is obtained by noting that
\begin{eqnarray*}
\triangledown w_{j} ( s ) & =&w_{j} ( s ) \triangledown\log
w_{j} ( s )
\\
& =&w_{j} ( s ) \bigl( v_{j}^{T}/c_{0}-w_{j}
( s ) v_{j}^{T}/c_{0} \bigr) =w_{j} (
s ) \bigl(1-w_{j} ( s ) \bigr)v_{j}^{T}/c_{0}.
\end{eqnarray*}
Therefore,
\[
\Delta\varrho_{b} ( s ) =\sum_{j=1}^{2m+d}w_{j}
( s ) \bigl( 1-w_{j} ( s ) \bigr) v_{j}v_{j}^{T}/c_{0}%
\]
and the result follows.
\end{pf}

Using the previous lemma it is a routine application of the dominated
convergence theorem to show that%
%
\begin{equation}
\triangledown H_{b}(s)=E \bigl( d^{\prime} \bigl(
\varrho_{b} ( s+X ) \bigr) \triangledown\varrho_{b} ( s+X )
\bigr) , \label{gradient}%
\end{equation}
where $d^{\prime} ( \cdot ) $ denotes the derivative of
$d (
\cdot ) $, and the gradient $\triangledown\varrho_{b} (
\cdot ) $ is encoded as a column vector. Similarly,
the Hessian matrix of $H_{b} ( \cdot ) $ is given by
%
\begin{equation}
\Delta H_{b}(s)=E \bigl( \Delta\varrho_{b} ( s+X )
d^{\prime
} \bigl( \varrho_{b} ( s+X ) \bigr) +d^{\prime\prime
}
\bigl( \varrho_{b} ( s+X ) \bigr) \triangledown\varrho _{b}
( s+X ) \triangledown\varrho_{b} ( s+X ) ^{T} \bigr) .
\label{hessian}%
\end{equation}
This will be useful in our technical development.

We now are ready to provide an estimate for the term $J_{2}$.

\begin{lemma}
\label{LemJ2} For every $\varepsilon\in(0,1/2),$ $\delta_{2}\in(0,1)$,
$c_{1}\in ( 1,\infty ) $ there exists $a\in(0,1)$ and $b_{0}>0$
(depending on $\varepsilon,\delta$, $v_{i}$'s, $\gamma$, $\delta
_{0}$) such
that if $b\geq b_{0}$, $s\notin b\Gamma\cup bA$ and $r_{b} (
s )
\leq-\delta_{2}b$ then
\[
J_{2} \bigl( 1-p_{b} ( s ) \bigr) \leq1-( 2-3\varepsilon )
\frac{P ( \exists j\dvt X^{T}v_{j}\geq a(a_{j}b-s^{T}v_{j}) ) } {H_{b}(s)} .
\]
\end{lemma}

\begin{pf}
Recall that%
%
\begin{eqnarray}
\label{loc} \bigl( 1-p_{b} ( s ) \bigr) J_{2} & =&E \biggl(
\frac
{g_{b}(s+X)}{g_{b}(s)};s+X\notin A_{b,a} ( s ) \biggr)
\nonumber
\\
& =&\int_{\mathbb{R}^{d}}\frac{g_{b}(s+x)}{g_{b}(s)}I\bigl(s+x\notin
A_{b,a} ( s ) \bigr)P ( X\in \mathrm{d}x )
\\
& =&\int_{\mathbb{R}^{d}}\frac{g_{b}(s+x)}{g_{b}(s)}I\bigl(s+x\notin
A_{b,a}(s), \Vert x\Vert _{2} <b \bigr)P ( X\in \mathrm{d}x ) +\mathrm{o}
\bigl(b^{-1}\bigr).\nonumber %
\end{eqnarray}
We will analyze the integrand above. Keep in mind that $s+x\notin
A_{b,a} ( s ) $. Given $c_{1}>0$ fixed, note that $g_{b}(s)<1$
whenever $r_{b}(s)\leq-\delta_{2}b$ assuming that $b$ is sufficiently large.
Therefore,
%
\begin{equation}
\frac{g_{b}(s+x)}{g_{b}(s)}\leq\frac{H_{b}(s+x)^{2}}{H_{b}(s)^{2}}%
=1+\frac{H_{b}(s+x)+H_{b}(s)}{H_{b}(s)}\times
\frac
{H_{b}(s+x)-H_{b}(s)}%
{H_{b}(s)}. \label{Ra}%
\end{equation}
In addition, applying a Taylor expansion, we obtain for each $x\in
\mathbb{R}^{d}$%
%
\begin{eqnarray}\label{HH}
H_{b}(s+x)-H_{b}(s) & =&\int_{0}^{1}x^{T}
\triangledown H_{b}%
(s+ux)\,\mathrm{d} u
\nonumber
\\[-8pt]
\\[-8pt]
\nonumber
& =&\int_{0}^{1}x^{T} \biggl(
\triangledown H_{b}(s)+u\int_{0}^{1}
\Delta H_{b}\bigl(s+u^{\prime}ux\bigr)x\,\mathrm{d} u^{\prime} \biggr)
\,\mathrm{d} u.
\end{eqnarray}
Observe that the previous expression can be written as%
%
\begin{equation}
H_{b}(s+x)-H_{b}(s)=x^{T}\triangledown
H_{b}(s)+E \bigl( x^{T}U\Delta H_{b}
\bigl(s+U^{\prime}Ux\bigr)x \bigr) , \label{SRep}%
\end{equation}
where $U$ and $U^{\prime}$ are i.i.d. $U ( 0,1 ) $. In what
follows, we consider the linear term $x^{T} \triangledown H_{b}(s)$ and the
quadratic term $x ^{T}u\Delta H_{b}(s+u^{\prime}ux)x$, respectively.

\textit{The linear term}.
We use the results in \eqref{gradient} and Lemma \ref{PropertiesRho}.
We shall
first start with the term involving $x^{T}\nabla H_{b} ( s ) $.
Define%
%
\begin{equation}
R=\bigl\{s\dvt r_{b} ( s ) \leq-\delta_{2}b,s\notin b
\Gamma\bigr\}. \label{R}%
\end{equation}
Observe that, uniformly over $s\in R$ we have%
\[
E \bigl( X^{T}\nabla H_{b} ( s ) I\bigl(s+X\notin
A_{b,a} ( s ) , \Vert x\Vert _{2} <b \bigr) \bigr) =\bigl(
\eta^{T}+\mathrm{o} \bigl( b^{-\alpha
+1}\bigr)\bigr)\triangledown
H_{b}(s).
\]
Furthermore,
\begin{eqnarray*}
\eta^{T}\triangledown H_{b}(s) & =&E \bigl(
\eta^{T}d^{\prime} \bigl( \varrho_{b}(s+X) \bigr)
\triangledown\varrho_{b}(s+X) \bigr)
\\
& =&E \Biggl( \sum_{j=1}^{2m+d}w_{j}(s+X)
\eta^{T}v_{j}^{\ast}d^{\prime
} \bigl(
\varrho_{b}(s+X) \bigr) \Biggr)
\\
& =&-E \bigl( d^{\prime} \bigl( \varrho_{b}(s+X) \bigr) \bigr).
\end{eqnarray*}
The last step is due to the fact that $\eta^{T}v_{j}=-1$ and $\sum_{j=1}^{2m+d}w_{j}(s+X)=1$. We have noted that%
%
\begin{eqnarray}\label{bd}
P \bigl( r_{b} ( s+X ) \geq\delta_{0} \bigr) & \leq& E
\bigl( d^{\prime} \bigl( \varrho_{b}(s+X) \bigr) \bigr)
\nonumber
\\[-8pt]
\\[-8pt]
\nonumber
& \leq& P \bigl( r_{b} ( s+X ) \geq-\delta_{0}-c_{0}
\log (2m+d) \bigr)
\nonumber
\end{eqnarray}
and therefore, if $s\in R$, following a reasoning similar to Lemmas
\ref{LemPropmu} and \ref{LemRatv}, for sufficiently large~$b_{0}$ (depending
on $\delta_{2},\gamma$) and $a$ sufficiently close to $1$ we have
that if
$b\geq b_{0}$
\[
a^{2\alpha}\leq\frac{P ( \exists j\dvt X^{T}v_{j}\geq a(a_{j}b-s^{T}
v_{j}) ) }{E ( d^{\prime} ( \varrho_{b}(s+X) )
)
}\leq a^{-2\alpha},
\]
which implies%
%
\begin{equation}
\eta^{T}\triangledown H_{b}(s)\leq-a^{2\alpha}P \bigl(
\exists j\dvt X^{T}%
v_{j}\geq a
\bigl(a_{j}b-s^{T}v_{j}\bigr) \bigr) .
\label{1stder}%
\end{equation}

Now, for $x$ taking values on any fixed compact set, we have uniformly over
$s\in R$ and $s+x\notin A_{b,a}(s)$ that%
%
\begin{equation}
\lim_{b\rightarrow\infty}\frac{H_{b}(s+x)+H_{b}(s)}{H_{b}(s)}=2; \label{lim2}%
\end{equation}
this follows easily using the bounds in (\ref{InR}) and the
representation in
Lemma \ref{LemRatv} in the \hyperref[SectionAppendixProperties]{Appendix}. We obtain that as $b\rightarrow
\infty
$,
\begin{eqnarray*}
&& \sup_{s\in R}\biggl\llvert E \biggl( \biggl( \frac{H_{b}(s+X)+H_{b}(s)} {H_{b}(s)}-2 \biggr) \times X^{T}
\frac{\triangledown
H_{b}(s)}{|\triangledown
H_{b}(s)|}I\bigl(s+X\notin A_{b,a} ( s ) , \Vert x\Vert
_{2} <b \bigr) \biggr) \biggr\rrvert
\\
&&\quad \leq E \biggl( \sup_{s\in R}\biggl\llvert \frac
{H_{b}(s+X)+H_{b}(s)}{H_{b}%
(s)}-2
\biggr\rrvert \times\frac{|X^{T}\triangledown
H_{b}(s)|}{|\triangledown
H_{b}(s)|}I\bigl(s+X\notin A_{b,a} ( s ) ,
\Vert x\Vert _{2} <b \bigr) \biggr)
\\
&&\quad =\mathrm{o} (1),
\end{eqnarray*}
where the last step is thanks to the uniform convergence in \eqref
{lim2} and
the dominated convergence theorem. Similar to the derivation of \eqref
{bd}, we
have that
\[
\bigl\llvert \triangledown H_{b}(s)\bigr\rrvert =\mathrm{O}  \bigl( P \bigl(
\exists j\dvt X^{T}v_{j}\geq a\bigl(a_{j}b-s^{T}v_{j}
\bigr) \bigr) \bigr) .
\]
Thus, we obtain that%

\begin{eqnarray}
\label{linear} && E \biggl( \frac{H_{b}(s+X)+H_{b}(s)}{H_{b}(s)}\times X^{T}\triangledown
H_{b}(s)I\bigl(s+X\notin A_{b,a} ( s ) , \Vert x\Vert
_{2} <b \bigr) \biggr)
\nonumber
\\[-8pt]
\\[-8pt]
\nonumber
&&\quad \leq\bigl(-2a^{2\alpha}+\mathrm{o} (1)\bigr)P \bigl( \exists j\dvt
X^{T}v_{j}\geq a\bigl(a_{j}%
b-s^{T}v_{j}\bigr) \bigr) , %
\end{eqnarray}
where the convergence corresponding to the term $\mathrm{o}  ( 1 ) $ is
uniform over $s\in R$ as $b\rightarrow\infty$.

\textit{The quadratic term}.
We proceed to consider $E ( x^{T}U\Delta H_{b}(s+U^{\prime
}Ux)x ) $
in (\ref{SRep}), which in turn feeds into (\ref{Ra}).


We take advantage of representation (\ref{SRep}). However, to avoid confusion
with the random variable $X$, we introduce an independent copy $\tilde
{X}$ of
$X$. Using \eqref{hessian}, we then obtain that for $u,u^{\prime}\in
(
0,1 ) $,%
\begin{eqnarray*}
&& \Delta H_{b}\bigl(s+u^{\prime}ux\bigr)
\\
&&\quad =E \bigl( \Delta\varrho_{b}\bigl(s+u^{\prime}ux+\tilde{X}
\bigr)d^{\prime
}\bigl(\varrho _{b}\bigl(s+u^{\prime}ux+
\tilde{X}\bigr)\bigr) \bigr)
\\
&&\qquad{}      +E \bigl( d^{\prime\prime}\bigl(\varrho_{b}
\bigl(s+u^{\prime}ux+\tilde {X}\bigr)\bigr)\triangledown
\varrho_{b}\bigl(s+u^{\prime}ux+\tilde {X}\bigr)\triangledown
\varrho_{b}\bigl(s+u^{\prime}ux+\tilde{X}\bigr)^{T}
\bigr)
\\
&&\quad =E \biggl( \frac{\sum_{j=1}^{2m+d}w_{j}(1-w_{j})v_{j}v_{j}^{T}}{c_{0}
}d^{\prime}\bigl(\varrho_{b}
\bigl(s+u^{\prime}ux+\tilde{X}\bigr)\bigr) \biggr)
\\
&&\qquad{}      +E \Biggl( \frac{I(|\varrho_{b}(s+u^{\prime}ux+\tilde
{X})|\leq
\delta_{0})}{2\delta_{0}}\sum_{j=1}^{2m+d}v_{j}w_{j}
\Biggl(\sum_{j=1} 
^{2m+d}v_{j}w_{j}
\Biggr)^{T} \Biggr) .
\end{eqnarray*}
Since $w_{j}\in(0,1)$, we have each of the element of $\Delta H_{b}%
(s+u^{\prime}ux)$ is bounded by%
%
\begin{eqnarray}
&& \frac{\sum_{j=1}^{2m+d}\Vert v_{j}\Vert _{2}^{2}}{c_{0}}P \bigl( \varrho _{b}\bigl(s+u^{\prime}ux+
\tilde{X}\bigr)\geq-\delta_{0} \bigr) \label {Term1}
\\
&&\quad{} +\frac{1}{2\delta_{0}} \Biggl(\sum_{j=1}^{2m+d}
\Vert v_{j}\Vert _{2} ^{2}%
\Biggr)^{2}P \bigl( \bigl|\varrho_{b}\bigl(s+u^{\prime}ux+
\tilde{X}\bigr)\bigr|\leq \delta _{0} \bigr) . \label{Term2}%
\end{eqnarray}

We will proceed to analyze each of the terms (\ref{Term1}) and (\ref{Term2})
separately. Let us start with (\ref{Term1}).\ Note that
%
\begin{equation}
\varrho_{b} \bigl( s+u^{\prime}ux+\tilde{X} \bigr) >-
\delta_{0} \label{IRhoJ2}%
\end{equation}
implies that
\[
r_{b} \bigl( s+u^{\prime}ux+\tilde{X} \bigr) \geq-
\delta_{0}-c_{0}%
\log(2m+d).
\]
Equivalently, there exists $j$ such that
\[
\tilde{X}^{T}v_{j}\geq a_{j}b-s^{T}v_{j}-u^{\prime}ux^{T}v_{j}-
\delta _{0}-c_{0}\log(2m+d).
\]
On the other hand, recall that $s+x\notin A_{b,a} ( s ) $ and
therefore, for every $j$%
\[
x^{T}v_{j}\leq a \bigl( a_{j}b-s^{T}v_{j}
\bigr) .
\]
Since $u,u^{\prime}\in(0,1)$, we have from our previous observations that
inequality (\ref{IRhoJ2}) implies that for some $j$
%
\begin{equation}
\tilde{X}^{T}v_{j}\geq(1-a) \bigl(a_{j}b-s^{T}v_{j}
\bigr)-\delta_{0}-c_{0}\log(2m+d). \label{CondtripStar}%
\end{equation}
Now, the inequalities%
\begin{eqnarray*}
r_{b} ( s ) & =&\max_{j=1}^{2m+d}
\bigl[s^{T}v_{j}-a_{j}b\bigr]
\\
& =&-\min_{j=1}^{2m+d}\bigl[a_{j}b-s^{T}v_{j}
\bigr]\leq-\delta_{2}b\leq -2\bigl(\delta _{0}+c_{0}
\log(2m+d)\bigr)/(1-a)
\end{eqnarray*}
together with (\ref{CondtripStar}) imply that for some $j$
%
\begin{equation}
\tilde{X}^{T}v_{j}\geq\tfrac{1}{2}(1-a)
\bigl(a_{j}b-s^{T}v_{j}\bigr)\geq
\tfrac
{1}%
{2}(1-a)\delta_{2}b. \label{Con1}%
\end{equation}
Therefore, we conclude that if $\delta_{2}b\geq2(\delta_{0}+c_{0}%
\log(2m+d))/(1-a)$, by regular variation there exists $c^{\prime}$
such that
for all $b$ sufficiently large%
\begin{eqnarray*}
&& P \bigl( \varrho_{b} \bigl( s+u^{\prime}ux+\tilde{X} \bigr)
>-\delta _{0} \bigr)
\\
&&\quad \leq P\bigl(\exists j\dvt \tilde{X}^{T}v_{j}\geq
\tfrac{1}{2}(1-a)\delta _{2}b\bigr)\leq c^{\prime}P\bigl(
\exists j\dvt \tilde{X}^{T}v_{j}\geq a\bigl(a_{j}b-s^{T}v_{j}
\bigr)\bigr)
\end{eqnarray*}
uniformly over $s\in R$. Note that in the last inequality we use the
fact that
if $s\notin b\Gamma\cup bA$ then $\Vert  s\Vert _{2}\leq cb$ for some constant $c$ depending on $\beta$ and
$\gamma$. Consequently, we conclude that (\ref{Term1}) is bounded by%
%
\begin{equation}
c^{\prime}\frac{\sum_{j=1}^{2m+d}\Vert v_{j}\Vert _{2} ^{2}}{c_{0}}P\bigl(\exists j\dvt X^{T}v_{j}^{\ast}
\geq a\bigl(a_{j}b-s^{T}v_{j}^{\ast}\bigr)
\bigr). \label{Term11} 
\end{equation}

Now we proceed with term (\ref{Term2}). Note that%
\[
P \bigl( \bigl\llvert \varrho_{b} \bigl( s+u^{\prime}ux+\tilde
{X} \bigr) \bigr\rrvert \leq\delta_{0} \bigr) \leq P \bigl(
\varrho_{b} \bigl( s+u^{\prime}ux+\tilde{X} \bigr) \geq-
\delta_{0} \bigr) ,
\]
so our previous development leading to the term (\ref{Term11})
implies that
given any selection of $\delta_{0}$, $a$, $c_{0}$ we have
%
\begin{equation}
\frac{P ( \llvert \varrho_{b} ( s+u^{\prime}ux+\tilde
{X} )
\rrvert \leq\delta_{0} ) }{P(\exists j\dvt X^{T}v_{j}\geq
a(a_{j}%
b-s^{T}v_{j}))}=\mathrm{O}  ( 1 ) \label{Term22a}%
\end{equation}
uniformly over $s+x\notin A_{b,a} ( s ) $ and $s\in R$ as
$b\rightarrow\infty$. Further, we have that as $b\rightarrow\infty$
\begin{eqnarray*}
\Delta_{s,b} \bigl( x,u,u^{\prime} \bigr) & \triangleq&
\frac{P (
\llvert \varrho_{b} ( s+u^{\prime}ux+\tilde{X} )
\rrvert
\leq\delta_{0} ) I ( s+x\notin A_{b,a} ( s )
)
}{P(\exists j\dvt X^{T}v_{j}\geq a(a_{j}b-s^{T}v_{j}))}
\\
& =&I \bigl( s+x\notin A_{b,a} ( s ) \bigr)
\\
&&{}      \times\frac{P ( \varrho_{b} ( s+u^{\prime
}ux+\tilde
{X} ) \geq\delta_{0} ) -P ( \varrho_{b} (
s+u^{\prime
}ux+\tilde{X} ) \geq-\delta_{0} ) }{P(\exists
j\dvt X^{T}v_{j}\geq
a(a_{j}b-s^{T}v_{j}))}
\\
& =&\mathrm{O} \bigl(b^{-1}\bigr).
\end{eqnarray*}
In addition, the above convergence is uniform on the set $s\in R$, and
$s+x\notin A_{b,a}(s)$, and on $u,u^{\prime}\in(0,1)$. We now
consider two
cases: $1<\alpha\leq2$ and $\alpha>2$.

\textit{Case}: $\alpha>2$.
The increment has finite second moment. Consider the quadratic term,
%
\begin{eqnarray}\label{quad}
&& \biggl\llvert \int_{\mathbb{R}^{d}}I \bigl( s+x\notin A_{b,a}
( s ) , \Vert x\Vert _{2} <b \bigr) \frac{H_{b}(s+x)+H_{b}(s)}{H_{b}(s)}E \bigl(
x^{T}U\Delta H_{b}\bigl(s+U^{\prime}Ux\bigr)x \bigr) P
( X\in \mathrm{d}x ) \biggr\rrvert
\nonumber\\
&&\quad \leq\frac{c^{\prime}}{c_{0}}E\Vert X\Vert _{2}^{2}\sum
_{j=1}^{2m+d}\Vert v_{j} 
\Vert _{2}^{2}P\bigl(\exists j\dvt X^{T}v_{j}
\geq a\bigl(a_{j}b-s^{T}v_{j}\bigr)\bigr)
\nonumber
\\
&&\qquad{} + P\bigl(\exists j\dvt X^{T}v_{j}\geq a
\bigl(a_{j}b-s^{T}v_{j}\bigr)\bigr)
\nonumber
\\[-8pt]
\\[-8pt]
\nonumber
&&\hphantom{\qquad{} +{} }{}     \times\int_{\mathbb{R}^{d}}I \bigl( s+x\notin A_{b,a} ( s
) , \Vert x\Vert _{2} <b \bigr) \frac{H_{b}(s+x)+H_{b}(s)}{H_{b}(s)}\\
&&  \hphantom{\qquad{} +{} {}     \times\int_{\mathbb{R}^{d}}}{}\times E\bigl\llvert
x^{T}\Delta_{s,b}\bigl(x,U,U^{\prime}\bigr)x\bigr\rrvert
P ( X\in \mathrm{d}x )
\nonumber
\\
&&\quad = \Biggl(\frac{c^{\prime}}{c_{0}}E\Vert X\Vert _{2}^{2}\sum
_{j=1}^{2m+d}\Vert v_{j}%
\Vert _{2}^{2}+\mathrm{o} (1) \Biggr)\times P\bigl(\exists j\dvt
X^{T}v_{j}\geq a\bigl(a_{j}b-s^{T}%
v_{j}\bigr)\bigr),
\nonumber
\end{eqnarray}
where the first step is thanks to the analysis results of \eqref
{Term1} and
\eqref{Term2} and the last step is thanks to the uniform convergence of
$\Delta_{s,b}(x,u,u^{\prime})$ and dominated convergence theorem.

\textit{Case}: $1< \alpha\leq2$.
The increment has infinite second moment. Note that for $s\in R$%
\[
\int_{\mathbb{R}^{d}}I \bigl( s+x\notin A_{b,a} ( s ) , \Vert x
\Vert _{2} <b \bigr) x^{T}x\,\mathrm{d} x=\mathrm{O} \bigl(b^{2-\alpha}\bigr).
\]
Given the choice of $c_{0}$ as in \eqref{c0}, the first term in the second
step of \eqref{quad} is
\begin{eqnarray*}
& &\frac{c^{\prime}}{c_{0}}\sum_{j=1}^{2m+d}\Vert
v_{j}\Vert _{2}^{2}P \bigl( \exists j\dvt
X^{T}v_{j}\geq a\bigl(a_{j}b-s^{T}v_{j}
\bigr) \bigr) \int_{\mathbb
{R}^{d}}I \bigl( s+x\notin A_{b,a} ( s
) , \Vert x\Vert _{2} <b \bigr) x^{T}x\,\mathrm{d} x
\\
&&\quad =\mathrm{o} (1)P\bigl(\exists j\dvt X^{T}v_{j}\geq a
\bigl(a_{j}b-s^{T}v_{j}\bigr)\bigr).
\end{eqnarray*}
In addition, given the convergence rate of $\Delta_{s,b}$, the second
term in
\eqref{quad} is
\begin{eqnarray*}
&& \int_{\mathbb{R}^{d}}I \bigl( s+x\notin A_{b,a} ( s ) , \Vert
x\Vert _{2} <b \bigr) \frac{H_{b}(s+x)+H_{b}(s)}{H_{b}(s)}E\bigl\llvert x^{T}
\Delta_{s,b}\bigl(x,U,U^{\prime}\bigr)x\bigr\rrvert P ( X\in \mathrm{d}x )
\\
&&\quad =\mathrm{O} \bigl(b^{1-\alpha}\bigr).
\end{eqnarray*}
Given that $c_{0}\geq b^{{(3-\alpha)}/{2}}$ and $b$ is large enough,
\eqref{quad} is bounded by
\[
\varepsilon P\bigl(\exists j\dvt X^{T}v_{j}\geq a
\bigl(a_{j}b-s^{T}v_{j}\bigr)\bigr)
\]
for all $\alpha>1$. Therefore, for all $\alpha>1$, \eqref{quad} is
bounded by
\[
2\varepsilon P\bigl(\exists j\dvt X^{T}v_{j}\geq a
\bigl(a_{j}b-s^{T}v_{j}\bigr)\bigr).
\]
%

\textit{Summary}.
We summarize the analysis of the linear and quadratic term by inserting bounds
in~\eqref{linear}, \eqref{quad}, and expansions \eqref{Ra} and \eqref
{SRep} into
\eqref{loc} and obtain that for $r_{b}(s)\leq-\delta_{2}b$ and
$s\notin
b\Gamma$
\begin{eqnarray*}
J_{2}\bigl(1-p_{b}(s)\bigr) & =&\int\frac{g_{b}(s+x)}{g_{b}(s)}I
\bigl( s+x\notin A_{b,a}(s) \bigr) P ( X\in \mathrm{d}x )
\\
& \leq&1-\bigl(2a^{2\alpha}+\mathrm{o} (1)+2\varepsilon\bigr)\frac{P(\exists
j\dvt X^{T}v_{j}\geq
a(a_{j}b-s^{T}v_{j}))}{H_{b}(s)}.
\end{eqnarray*}
We conclude the proof of the lemma by choosing $a$ sufficiently close
to $1$
and $b$ sufficiently large.
\end{pf}

The subsequent lemmas provide convenient bounds for $p_{b} (
s ) $,
$J_{1}$ and $J_{2}$.

\begin{lemma}
\label{Lempb} For each $\varepsilon>0$ and any given selection of
$\tilde
{c}_{0}$ (feeding into $c_{0}$), $\delta_{0},\gamma,\delta_{2}>0$
there exists
$b_{0}\geq0$ such that if $b\geq b_{0}$, $r_{b} ( s ) \leq
-\delta_{2}b$, and $s\notin b\Gamma$ then
\begin{longlist}[(ii)]
\item[(i)]
\[
p_{b} ( s ) \leq\varepsilon.
\]

\item[(ii)]
\[
1\leq\frac{H_{b}(s)}{v_{b} ( s ) }\leq1+\varepsilon.
\]
\end{longlist}
\end{lemma}

\begin{pf}
Note that if $r_{b} ( s ) =\max_{j=1}^{2m+d}[s^{T}v_{j}-a_{j}%
b]\leq-\delta_{2}b$ and there exists $j$ such that $X^{T}v_{j}\geq
a(a_{j}b-s^{T}v_{j})$, we must also have that there exists $j$ such that
\[
X^{T}v_{j}\geq a\bigl(a_{j}b-s^{T}v_{j}
\bigr)\geq a\delta_{2}b.
\]
On the other hand, $s\notin b\Gamma\cup bA$ then $\Vert
s\Vert _{2}=\mathrm{O}  ( b ) $, together with
$c_{0} =
\mathrm{o} (b)$, there is a constant $c$ (depending on $\beta$ and $\gamma
$) such
that%
\begin{eqnarray*}
H_{b}(s) & \geq& E\bigl[\varrho_{b}^{+} ( s+X )
\bigr]=\int_{0}^{\infty
}P \bigl( \varrho_{b}
( s+X ) >t \bigr) \,\mathrm{d} t
\\
& \geq&\int_{0}^{\infty}P \bigl( r_{b} (
s+X ) >t - c_{0} \log(2m+d) \bigr) \,\mathrm{d} t
\\
& =&\int_{0}^{\infty}P \bigl( \exists j\dvt
X^{T}v_{j}\geq a\bigl(a_{j}b-s^{T}v_{j}
\bigr)+t- c_{0} \log(2m+d) \bigr) \,\mathrm{d} t
\\
& \geq&\int_{0}^{\infty}P \bigl( \exists j\dvt
X^{T}v_{j}\geq cb+t \bigr) \,\mathrm{d} t
\\
& =&b\int_{0}^{\infty}P \bigl( \exists j\dvt
X^{T}v_{j}\geq cb+ub \bigr) \,\mathrm{d} u
\geq\delta^{\prime}bP \bigl( \exists j\dvt X^{T}v_{j}
\geq(c+1)b \bigr)
\end{eqnarray*}
for some $\delta^{\prime}>0$ small. Therefore,
\[
p_{b} ( s ) \leq\frac{\theta P ( \exists
j\dvt X^{T}v_{j}\geq
a\delta_{2}b ) }{\delta^{\prime}bP ( \exists j\dvt X^{T}v_{j}%
\geq(c+1)b ) }=\mathrm{O}  ( 1/b )
\]
as $b\rightarrow\infty$; this bound immediately yields (i). Part (ii) is
straightforward from the estimates in~(\ref{SInrho}), (\ref{InR}), and
eventually \eqref{UB}.
\end{pf}

\begin{lemma}
\label{LemJ1} For any $\varepsilon>0$, we can choose $b$ sufficiently
large so
that on the set $r_{b}(s)\leq-\delta_{2}b$,
%
\[
J_{1}\leq(1+\varepsilon)\frac{P ( \exists j\dvt X^{T}v_{j}\geq
a(a_{j}%
b-s^{T}v_{j}) ) }{c_{1}\theta H_{b}(s)}.
\]
\end{lemma}

\begin{pf}
Choose $b$ large enough such that $g_{b}(s)<1$ whenever $r_{b}(s)\leq
-\delta_{2}b$. With $p_{b}(s)$ defined as in \eqref{D-p} and the fact that
$g_{b}(s+X)\leq1$, it follows easily that
\begin{eqnarray*}
J_{1} & \leq&\frac{P ( \exists j\dvt X^{T}v_{j}\geq a(a_{j}b-s^{T}%
v_{j}) ) ^{2}}{c_{1}H_{b}(s)^{2} ( p_{b} ( s )
+ (
1-p_{b} ( s )  ) P ( \exists j\dvt X^{T}v_{j}\geq
a(a_{j}b-s^{T}v_{j}) )  ) }
\\
& \leq&\frac{P ( \exists j\dvt X^{T}v_{j}\geq
a(a_{j}b-s^{T}v_{j}) )
^{2}}{c_{1}H_{b}(s)^{2}p_{b} ( s ) }
 =(1+\varepsilon) \frac{P ( \exists j\dvt X^{T}v_{j}\geq a(a_{j}%
b-s^{T}v_{j}) ) }{c_{1}\theta H_{b}(s)}.
\end{eqnarray*}
The last step uses the approximation \eqref{UB} and the definition of
$p_{b}(s)$ in \eqref{D-p}.
\end{pf}

We summarize our verification of the validity of $g_{b}(s) $ in the
next result.

\begin{proposition}
\label{Prop}Given $\varepsilon,\delta_{2}>0$, let $a$ be as chosen
in Lemma
\ref{LemJ2}. We choose $\theta=1/(1+\varepsilon)^{2}$ and
$c_{1}= (
1+\varepsilon ) ^{3}(1+4\varepsilon)$ such that on the set $r_{b}(s)
\leq-\delta_{2}b$ and $s\notin b\Gamma$ we have
\[
J_{1}+J_{2}\leq1.
\]
\end{proposition}

\begin{pf}
Combining our bounds for $p_{b} ( s ) $, $J_{1}$ and $J_{2}$ given
in Lemmas \ref{LemJ2}, \ref{Lempb}, and \ref{LemJ1}, we obtain
\begin{eqnarray*}
&& J_{1}+J_{2}
\\[-2pt]
&&\quad \leq1+p_{b} ( s ) +\varepsilon p_{b} ( s ) +(1+
\varepsilon)\frac{P ( \exists j\dvt X^{T}v_{j}\geq a(a_{j}b-s^{T}
v_{j}) ) }{c_{1}\theta H_{b}(s)}
\\[-2pt]
&&\qquad{} -2\frac{P ( \exists j\dvt X^{T}v_{j}\geq a(a_{j}b-s^{T}v_{j}) )
}{H_{b} ( s ) } ( 1-3\varepsilon/2 )
\\[-2pt]
& &\quad\leq1+\frac{P ( \exists j\dvt X^{T}v_{j}\geq
a(a_{j}b-s^{T}v_{j}) )
}{H_{b} ( s ) } \biggl( \theta(1+\varepsilon)^{2}+
\frac
{1+\varepsilon}{c_{1}\theta} \biggr)
\\[-2pt]
&&\qquad{} -2\frac{P ( \exists j\dvt X^{T}v_{j}\geq a(a_{j}b-s^{T}v_{j}) )
}{H_{b} ( s ) } ( 1-3\varepsilon/2 ) .
\end{eqnarray*}
We arrive at
\begin{eqnarray*}
&& J_{1}+J_{2}
\\[-2pt]
&&\quad \leq1+\frac{P ( \exists j\dvt X^{T}v_{j}\geq
a(a_{j}b-s^{T}v_{j}) )
}{H_{b} ( s ) } \biggl( \theta(1+\varepsilon)^{2}+
\frac
{1+\varepsilon}{c_{1}\theta}-2 ( 1-3\varepsilon/2 ) \biggr) .
\end{eqnarray*}
We then select $\theta=1/(1+\varepsilon)^{2}$ and $c_{1}= (
1+\varepsilon ) ^{3}(1+4\varepsilon)$ and conclude that
$J_{1}+J_{2}%
\leq1$ as required for if $\varepsilon\in(0,1/2)$.
\end{pf}

\subsection{\texorpdfstring{Proofs of Theorem \protect\ref{ThmMain1} and Proposition \protect\ref{PropFinieH}}
{Proofs of Theorem 1 and Proposition 1}}

We now are ready to provide the proof of Proposition \ref
{PropFinieH}. As
noted earlier, Theorem \ref{ThmMain1} follows immediately as a
consequence of
Proposition \ref{PropFinieH}, combined with Lemma \ref{LemFiniteHorizon}.

\begin{pf*}{Proof of Proposition \ref{PropFinieH}}
For each $\lambda\in (
-\infty,\infty ) $, define $\tau ( \lambda ) =\inf
\{n\geq0\dvt r_{b} ( S_{n} ) \geq\lambda\}$. Proposition \ref{Prop}
together with Lemma \ref{LemLB} implies
\[
g_{b} ( s ) \geq E_{s} \Biggl( g_{b} (
S_{\tau (
-\delta_{2}b ) } ) \prod_{j=0}^{\tau ( -\delta
_{2}b ) -1}\hat{k} (
S_{j},S_{j+1} ) I \bigl( \tau ( -\delta_{2}b )
\leq T_{b\Gamma} \bigr) \Biggr) .
\]
Given $c_{1}>0$, there exists $\kappa\in ( 0,\infty ) $
such that
$r_{b} ( s ) \geq\kappa$ implies $g ( s ) =1$.
So, we
have that%
%
\begin{eqnarray}\label{Bnd1}
g_{b} ( s ) & \geq& E_{s} \Biggl( g_{b} (
S_{\tau (
-\delta_{2}b ) } ) \prod_{j=0}^{\tau ( -\delta
_{2}b ) -1}\hat{k} (
S_{j},S_{j+1} ) I \bigl( \tau ( -\delta_{2}b )
\leq T_{b\Gamma} \bigr) \Biggr)
\nonumber
\\[-8pt]
\\[-8pt]
\nonumber
& \geq &E_{s} \Biggl( \prod_{j=0}^{\tau ( -\delta_{2}b )
-1}
\hat{k} ( S_{j},S_{j+1} ) I \bigl( \tau ( -\delta
_{2}b ) \leq T_{b\Gamma},r_{b} ( S_{\tau ( -\delta
_{2}b ) } )
\geq\kappa \bigr) \Biggr) .
\nonumber
\end{eqnarray}
We prove the theorem in two steps. The first step is to show that $\hat
{P}%
_{0} ( \cdot ) $ is a good approximation, in total
variation, to
$P_{0}(\cdot|\tau ( -\delta_{2}b ) <T_{b\Gamma},r_{b}(
S_{\tau(
-\delta_{2}b) }) \geq\kappa)$. In step 2, we show that $P_{0}(\cdot
|\tau ( -\delta_{2}b ) <T_{b\Gamma},r_{b} ( S_{\tau
(
-\delta_{2}b ) } ) \geq\kappa)$ approximates
$P_{0}(\cdot
|T_{bA}<T_{b\Gamma})$ well.

\textit{Step} 1.
Applying Lemmas \ref{LemTV1}, \ref{LemLB}, and a similar
argument to the
one below Lemma \ref{LemTV1}, it is sufficient to show that for any
$\varepsilon>0$ we can pick $\delta_{2}$, $\gamma$, $\theta$, $a$, $c_{1}$
such that for all $b$ large enough%
\[
g_{b}(0)\leq(1+\varepsilon)P_{0}^{2}\bigl(
\tau(-\delta_{2}b)<T_{b\Gamma} ,r_{b}(S_{\tau(-\delta_{2}b)})
\geq\kappa\bigr).
\]
First, note that%
\[
P_{0} \bigl( \tau ( -\delta_{2}b ) <\infty \bigr) \geq
H_{b}(0) .
\]
Therefore,%
\begin{eqnarray*}
& &P_{0} \bigl( \tau ( -\delta_{2}b ) <\infty
,r_{b}(S_{\tau
(-\delta_{2}b)})\geq\kappa \bigr)
\\
&&\quad            =P_{0} \bigl( \tau ( -\delta_{2}b ) <\infty
\bigr) \times P_{0} \bigl( r_{b}(S_{\tau(-\delta_{2}b)})\geq
\kappa|\tau ( -\delta_{2}b ) <\infty \bigr)
\\
&&\quad            \geq H_{b}(0)\times P_{0} \bigl(
r_{b} ( S_{\tau
(
-\delta_{2}b ) } ) \geq\kappa|\tau ( -\delta _{2}b
) <\infty \bigr) .
\end{eqnarray*}
Then, given $c_{1}=(1+\varepsilon)$, we show that one can pick $\delta_{2}>0$
small enough, depending on $\varepsilon>0$, such that for $b$ sufficiently
large
%
\begin{equation}
P_{0} \bigl( r_{b}(S_{\tau(-\delta_{2}b)})\geq\kappa|\tau ( -
\delta _{2}b ) <\infty \bigr) \geq1/(1+\varepsilon).
\label{Clmd2}%
\end{equation}
In order to do this we will use results from one dimensional regularly varying
random walks. Define%
\[
\tau_{i} ( -\delta_{2}b ) =\inf\bigl\{n\geq 0\dvt
v_{i}^{T}S_{n}-a_{i}%
b\geq-
\delta_{2}b\bigr\}
\]
for $i=1,\ldots,2m+d$ and let
\[
\mathcal{I}=\Bigl\{1\leq i\leq2m+d\dvt \mathop{\overline{\lim}}_{b\rightarrow\infty
}P (
T_{bA}<\infty ) /P \bigl( \tau_{i}(0)<\infty \bigr) <\infty
\Bigr\}.
\]
Observe that%
\[
\max_{i=1}^{2m+d}P \bigl( \tau_{i}(0)<
\infty \bigr) \leq P ( T_{bA}<\infty ) \leq(2m+d)\max
_{i=1}^{2m+d}P \bigl( \tau_{i}%
(0)<\infty \bigr) ,
\]
so the set $\mathcal{I}$ contains the half-spaces that have substantial
probability of being reached given that $T_{bA}<\infty$. We have that
%
\begin{eqnarray}
\label{B1dim} && \frac{P_{0} ( r_{b} ( S_{\tau ( -\delta_{2}b )
} ) \leq\kappa,\tau ( -\delta_{2}b ) <\infty
) }%
{P_{0}\bigl(\tau ( -
\delta_{2}b ) <\infty\bigr)}
\nonumber
\\
&&\quad \leq\sum_{i=1}^{2m+d}\frac{P_{0} ( r_{b} ( S_{\tau (
-\delta_{2}b ) } ) \leq\kappa,\tau_{i}(-\delta
_{2}b)=\tau (
-\delta_{2}b ) ,\tau ( -\delta_{2}b ) <\infty
) }%
{P_{0}\bigl(\tau ( -\delta_{2}b ) <\infty\bigr)}
\nonumber
\\[-8pt]
\\[-8pt]
\nonumber
&&\quad \leq\sum_{i=1}^{2m+d}\frac{P_{0} ( r_{b} ( S_{\tau
_{i} (
-\delta_{2}b ) } ) \leq\kappa,\tau_{i} ( -\delta
_{2}b )
<\infty ) }{P_{0}(\tau ( -\delta_{2}b ) <\infty
)}
\nonumber\\
&&\quad \leq\sum_{i\in\mathcal{I}}\frac{P_{0} ( v_{i}^{T}S_{\tau
_{i} (
-\delta_{2}b ) }-a_{i}b\leq\kappa,\tau_{i} ( -\delta
_{2}b )
<\infty ) }{P_{0}(\tau ( -\delta_{2}b ) <\infty
)}+\mathrm{o}  ( 1 ) .\nonumber
\end{eqnarray}
Now, $i\in\mathcal{I}$ implies that $v_{i}^{T}X$ is regularly varying with
index $\alpha$ and therefore (see \cite{AsmKlup96}) there exists a constant
$c_{i}^{\prime}>0$ (independent of $\delta_{2}$) such that for all
$b$ large
enough%
\[
P_{0} \bigl( v_{i}^{T}S_{\tau_{i} ( -\delta_{2}b ) }-(a_{i}
-\delta_{2})b>2\delta_{2}b|\tau_{i}
( -\delta_{2}b ) <\infty \bigr) \geq\frac{1}{ ( 1+2c_{i}^{\prime}\delta_{2} )
^{\alpha-1}%
}.
\]
Consequently,
%
\begin{eqnarray}
\label{Ind1stpas} && \frac{P_{0} ( v_{i}^{T}S_{\tau_{i} ( -\delta_{2}b
) }%
-a_{i}b\leq\kappa,\tau_{i} ( -\delta_{2}b ) <\infty
) }%
{P_{0}\bigl(\tau ( -
\delta_{2}b ) <\infty\bigr)}
\nonumber
\\[-8pt]
\\[-8pt]
\nonumber
&&\quad \leq\bigl(1- \bigl( 1+2c_{i}^{\prime}\delta_{2}
\bigr) ^{1-\alpha
}\bigr)\frac
{P_{0} ( \tau_{i} ( -\delta_{2}b ) <\infty ) }%
{P_{0} \bigl(
\tau ( -\delta_{2}b ) <\infty \bigr) }. %
\end{eqnarray}
Bound of (\ref{Ind1stpas}) together with (\ref{B1dim}) implies
that one can
select $\delta_{2}$ small enough and $b$ large enough to satisfy
(\ref{Clmd2}) and thus
\[
P_{0} \bigl( \tau ( -\delta_{2}b ) <\infty,r_{b}(S_{\tau
(-\delta_{2}b)})
\geq\kappa \bigr) \geq(1- \varepsilon)H_{b}(0).
\]
Furthermore, one can choose $\gamma$ large enough so that
\[
P_{0} \bigl( \tau ( -\delta_{2}b ) <T_{b\Gamma},r_{b}%
(S_{\tau(-\delta_{2}b)})\geq\kappa \bigr) \geq(1- \varepsilon)H_{b}(0).
\]
Thanks to Lemma \ref{LemTV1}, we obtain that
%
\begin{equation}
\label{s1}\mathop{\overline{\lim}}_{b\rightarrow\infty}\sup_{B}\bigl\llvert
\hat P_{0}(B) -P_{0}\bigl( B|\tau(-\delta_{2}b)
< T_{b\Gamma},r_{b}( S_{\tau( -\delta_{2}b)}) \geq\kappa\bigr) \bigr
\rrvert \leq\varepsilon.
\end{equation}
This completes step 1.

\textit{Step} 2.
With an analogous development for \eqref{Clmd2}, we can establish
that with
$\delta_{2}$ chosen small enough and $\gamma$ large enough
%
\begin{equation}
\label{s2}\mathop{\overline{\lim}}_{b\rightarrow\infty}\sup_{B}\bigl\llvert
P_{0} \bigl( B|\tau ( -\delta_{2}b ) <T_{b\Gamma}
\bigr) -P_{0} \bigl( B|\tau ( -\delta_{2}b ) <
T_{b\Gamma},r_{b} ( S_{\tau (
-\delta_{2}b ) } ) \geq\kappa \bigr) \bigr
\rrvert \leq \varepsilon.
\end{equation}
In addition, $\{T_{bA}<T_{b\Gamma}\}\subset\{\tau ( -\delta
_{2}b )
<T_{b\Gamma}\}$. Using the same type of reasoning leading to
(\ref{Ind1stpas}) we have that for any $\varepsilon>0$, $\delta
_{2}$ can be
chosen so that
\[
P_{0}(T_{bA}<T_{b\Gamma})\geq\frac{1}{1+\varepsilon}P_{0}
\bigl(\tau (-\delta _{2}b)<T_{b\Gamma}\bigr)
\]
for $b$ and $\gamma$ large enough. Therefore, we have
%
\begin{equation}
\mathop{\overline{\lim}}_{b\rightarrow\infty}\sup_{B}\bigl\llvert
P_{0} \bigl( B|\tau( -\delta_{2}b) <T_{b\Gamma} \bigr)
-P_{0}(B|T_{bA}<T_{b\Gamma}) \bigr\rrvert \leq
\varepsilon. \label{3}%
\end{equation}
Combining \eqref{s2} and \eqref{3}, we obtain that with $\gamma$
large enough
%
\begin{equation}
\mathop{\overline{\lim}}_{b\rightarrow\infty}\sup_{B}\bigl\llvert
P_{0} \bigl( B|\tau ( -\delta_{2}b ) < T_{b\Gamma},r_{b}
( S_{\tau (
-\delta_{2}b ) } ) \geq\kappa \bigr) -P_{0}(B|T_{bA}<T_{b\Gamma})
\bigr\rrvert \leq\varepsilon. \label{s3}%
\end{equation}

Then, we put together \eqref{s1}, \eqref{s3} and conclude that
\[
\mathop{\overline{\lim}}_{b\rightarrow\infty}\sup_{B}\bigl\llvert \hat
{P}_{0}%
(B)-P_{0} ( B|T_{bA}<T_{b\Gamma}
) \bigr\rrvert \leq 2\varepsilon.
\]
\upqed\end{pf*}

\section{A conditional central limit theorem}\label{SectCCLT}

The goal of this section is to provide a proof of Theorem \ref
{ThmMain2}. The
proof is a consequence of the following proposition. First, write for any
$t\geq0$%
%
\begin{equation}
\kappa_{\mathbf{a}}(t)\triangleq\mu \Bigl( \Bigl\{y\dvt \max
_{j=1}^{2m+d} \bigl( y^{T}v_{j}-a_{j}
\bigr) >t\Bigr\} \Bigr) , \label{DefKa}%
\end{equation}
where $\mathbf{a}=(a_{1},\ldots,a_{2m+d})$. Recall that $A$ is defined as
in \eqref{A}.

\begin{proposition}
\label{PropCLTAux}For each $z>0$ let $Y ( z ) $ be a random
variable with distribution given by
\[
P \bigl( Y ( z ) \in B \bigr) =\frac{\mu ( B\cap
\{y\dvt \max_{j=1}^{2m+d}[y^{T}v_{j}-a_{j}]\geq z\} ) }{\mu (
\{y\dvt \max_{j=1}^{2m+d}[y^{T}v_{j}-a_{j}]\geq z\} ) }.
\]
In addition, let $Z$ be a positive random variable following distribution
\[
P ( Z>t ) =\exp \biggl\{ -\int_{0}^{t}
\frac{\kappa
_{\mathbf{a}}%
(s)}{\int_{s}^{\infty}\kappa_{\mathbf{a}}(u)\,\mathrm{d} u}\,\mathrm{d} s \biggr\} .
\]
Let $\hat S_{0}=0$ and $\hat S_{n}$ evolve according to the transition kernel
\eqref{DKt} associated with the approximation in Proposition
\ref{PropFinieH}. Define $T_{bA} = \inf\{n\dvt \hat S_{n} \in bA\}$.
Then, as
$b\rightarrow\infty$, we have that
\[
\biggl( \frac{T_{bA}}{b},\frac{\hat S_{uT_{bA}}-u T_{bA}\eta}{\sqrt {T_{bA}}%
},\frac{\hat X_{T_{bA}}}{b} \biggr)
\Rightarrow \bigl( Z,CB ( uZ ) ,Y ( Z ) \bigr)
\]
in $\mathbb{R}\times D[0,1)\times\mathbb{R}^{d}$, where $CC^{T}=\operatorname{Var}(X)$,
$B ( \cdot ) $ is a $d$-dimensional Brownian motion with identity
covariance matrix, $B(\cdot)$ is independent of $Z$ and $Y (
Z ) $.
\end{proposition}

The strategy to prove Proposition \ref{PropCLTAux} is to create a coupling
of two processes $S$ and $\hat{S}$ on the same probability space with measure
$P$. For simplicity, we shall assume that $\hat{S}_{0}=S_{0}=0$. The process
$S=\{S_{n}\dvt n\geq0\}$ follows its original law, that is, $S_{n}=X_{1}%
+\cdots+X_{n}$ where $X_{i}$'s are i.i.d. The process $\hat{S}$ evolves
according to the transition kernel
\[
P(\hat{S}_{n+1}\in \mathrm{d}s_{n+1}|\hat{S}_{n}=s_{n})=
\hat{K}(s_{n},\mathrm{d}s_{n+1})
\]
obtained in Theorem \ref{ThmMain1}. Now we explain how the process $S$ and
$\hat{S}$ are coupled. A transition at time $j$, given $\hat
{S}_{j-1}=\hat
{s}_{j-1}$, is constructed as follows. First, we construct a Bernoulli random
variable $I_{j}$ with success parameter $p_{b} ( \hat
{s}_{j-1} ) $.
If $I_{j}=1$ then we consider $X$ (a generic random variable following the
nominal/original distribution) given that $\hat{s}_{j-1}+X\in A_{b,a}
(\hat{s}_{j-1})$ and let $\hat{X}_{j}=X$ (recall that $A_{b,a}(s)$ is defined
as in \eqref{Aa}); otherwise if $I_{j}=0$ we let $\hat{X}_{j}=X_{j}$,
that is,
we let $\hat{X}_{j}$ be equal to the $j$th increment of $S$. We then define
$N_{b}=\inf\{n\geq1\dvt I_{n}=1\}$ and observe that $\hat{S}_{j}=S_{j}$ for
$j<N_{b}$. The increments of processes $S$ and $\hat{S}$ at times $j>N_{b}$
are independent.

We will first show that $N_{b}=T_{bA}$ with high probability as
$b\nearrow
\infty$, once this result has been shown the rest of the argument basically
follows from functional central limit theorem for standard random walk. We
need to start by arguing that whenever a jump occurs (i.e., $I_{n}=1$)
the walk
reaches the target set with high probability. This is the purpose of the
following lemma.

\begin{lemma}
\label{LemJump}For every $\varepsilon$, $\delta_{2}$, $\gamma>0$
there exists
$a,b_{0}>0$ such%
\[
P \bigl( r_{b} ( s+X ) >0\vert s+X\in A_{b,a} ( s ) \bigr)
\geq1-\varepsilon
\]
for all $r_{b} ( s ) \leq-\delta_{2}b$, $s\notin b\Gamma
$, and $b>
b_{0}$.
\end{lemma}

\begin{pf}
Set $s=b\cdot u\in\mathbb{R}^{d}$ and note that%
\begin{eqnarray*}
&& P \Bigl( \max_{j=1}^{2m+d} \bigl( ( s+X )
^{T}%
v_{j}-a_{j}b \bigr) >0\vert s+X
\in A_{b,a} ( s ) \Bigr)
\\
&&\quad =P \bigl( \exists j \dvt X^{T}v_{j}\geq b
\bigl(a_{j}%
-u^{T}v_{j}\bigr)\vert
\exists j \dvt X^{T}v_{j}\geq ab \bigl( a_{j}-u^{T}v_{j}
\bigr) \bigr)
\\
&&\quad =\frac{P ( \exists j\dvt X^{T}v_{j}\geq
b(a_{j}-u^{T}%
v_{j}) ) }{P ( \exists j\dvt X^{T}v_{j}\geq
ab (
a_{j}-u^{T}v_{j} )  ) }.
\end{eqnarray*}
For each fixed $u,$ we have that%
\[
\lim_{b\rightarrow\infty}\frac{P ( \exists j \dvt X^{T}%
v_{j}\geq b(a_{j}-u^{T}v_{j}) ) }{P ( \exists j\dvt X^{T}v_{j}\geq ab ( a_{j}-u^{T}v_{j} )  ) }=\frac
{\mu\{y\dvt \exists j\dvt y^{T}v_{j}\geq
(a_{j}-u^{T}v_{j})\}}%
{
\mu\{y\dvt \exists j \dvt y^{T}v_{j}\geq a
(a_{j}-u^{T}v_{j})\}}.
\]
The convergence occurs uniformly over the set of $u$'s such that
$r_{b} (
ub ) /b\leq-\delta_{2}$ and $u\notin\Gamma$, which is a
compact set. The
result then follows by continuity of the radial component in the polar
representation of $\mu ( \cdot )$ (Lemma \ref
{LemPropmu}) as
$a\rightarrow1$.
\end{pf}

Now we prove that $T_{bA}=N_{b}$ occurs with high probability as
$b\nearrow\infty$.

\begin{lemma}
\label{LemTeqN1}For any $\varepsilon>0,$ we can select $\gamma>0$
sufficiently
large so that
\[
\mathop{\underline{\lim}}_{b\rightarrow\infty}P_{0}(N_{b}<\infty)\geq 1-
\varepsilon.
\]
\end{lemma}

\begin{pf}
Define $T_{b\Gamma}=\inf\{n\dvt \hat S_{n}\in b\Gamma\}$. Choose
$\varepsilon
^{\prime}$ positive and $\gamma$ large enough so that for all $t <
(1-\varepsilon^{\prime})\gamma/d$
\begin{eqnarray*}
P_{0} ( N_{b}>tb ) & \leq& P(N_{b} > tb,
T_{b\Gamma}> tb) + \mathrm{o} (1)
\\
& =&E_{0} \biggl[ \prod_{k\leq tb}
\bigl(1-p_{b}(\hat{S}_{k})\bigr);T_{b\Gamma}> tb
\biggr] +\mathrm{o} (1)
\\
& \leq& E_{0} \biggl[ \prod_{k\leq tb}
\bigl(1-p_{b}(\hat{S}_{k})\bigr)I\bigl(\Vert \hat
{S}_{k}-\eta k\Vert _{2}<\varepsilon^{\prime}\max
\{k,b\}\bigr) \biggr]
\\
&&{} +P_{0} \Bigl( \sup_{k\leq tb}\bigl\llvert \Vert
\hat{S}_{k}-\eta k\Vert _{2}%
-
\varepsilon^{\prime}\max\{k,b\}\bigr\rrvert >0 \Bigr) +\mathrm{o} (1).
\end{eqnarray*}
In the last step, we drop the condition $T_{b\Gamma}> tb$ on the set
$\{\Vert \hat{S}_{k}-\eta k\Vert _{2}<\varepsilon^{\prime}\max\{k,b\}\}$.
The second
term in the last step vanishes as $b\rightarrow\infty$ for any
$\varepsilon
^{\prime}>0$.

We claim that we can find constants $\delta^{\prime},c^{\prime}>0$
such that
for all $1\leq k\leq c^{\prime}\gamma b$%
%
\begin{equation}
\inf_{\{s\dvt \Vert \eta k-s\Vert _{2}\leq\varepsilon
^{\prime}%
\max\{k,b\},s\notin b\Gamma\}}p_{b}(s)\geq\frac{\delta^{\prime}}{k+b}.
\label{Claimaa}%
\end{equation}
To see this, recall that if $r_{b} ( s ) \leq-\delta_{2}b$
then%
\[
p_{b}(s)=\frac{\theta P ( \exists j\dvt X^{T}v_{j}>a(a_{j}b-s^{T}%
v_{j}) ) }{\int_{0}^{\infty}P ( \exists j\dvt X^{T}v_{j}>a_{j}%
b-s^{T}v_{j}+t ) \,\mathrm{d} t}.
\]
Now, if $\Vert  \eta k-s\Vert _{2}\leq\varepsilon^{\prime}\max\{k,b\}$
then, letting
$\lambda_{+}=\max_{j=1}^{2m+d}\Vert  v_{j}\Vert _{2}$ we obtain, by the Cauchy--Schwarz inequality%
%
\begin{eqnarray}
\bigl\llvert s^{T}\eta-kd\bigr\rrvert & \leq& d^{1/2}
\varepsilon ^{\prime}%
\max\{k,b\},\label{InAuxa}
\\
\bigl\llvert s^{T}v_{j}-k\eta^{T}v_{j}
\bigr\rrvert & \leq&\lambda \varepsilon ^{\prime}\max\{k,b\}.
\label{InAuxb}%
\end{eqnarray}
Inequality (\ref{InAuxa}) implies that
\[
s^{T}\eta\leq d^{1/2}\varepsilon^{\prime}\max\{k,b\}+kd
\leq k\bigl(\varepsilon ^{\prime} d ^{1/2}+d\bigr)+
\varepsilon^{\prime}d ^{1/2}b.
\]
We choose $\varepsilon^{\prime}d ^{1/2}<\gamma/2$. Then $k\leq
\gamma
b/2(\varepsilon^{\prime}d ^{1/2}+d)$ implies $s^{T}\eta<\gamma b$.
We shall
select
\[
c^{\prime}=\frac1{2\bigl(\varepsilon^{\prime}d ^{1/2}+d
\bigr)}.
\]
Inequality (\ref{InAuxb}) implies that
\begin{eqnarray*}
p_{b}(s) & \geq&\frac{\theta P ( \exists j\dvt X^{T}v_{j}>a(a_{j}%
b+k+\lambda\varepsilon^{\prime}\max\{k,b\}) ) }{\int_{0}^{\infty
}P ( \exists j\dvt X^{T}v_{j}>a_{j}b+k+t-\lambda\varepsilon^{\prime}
\max\{k,b\} ) \,\mathrm{d} t}
\\
& \geq&\frac{\theta P ( \exists j\dvt X^{T}v_{j}>a((a_{j}+\lambda
\varepsilon^{\prime})b+k(1+\lambda\varepsilon^{\prime})) )
}{\int_{0}^{\infty}P ( \exists j\dvt X^{T}v_{j}>(a_{j}-\lambda\varepsilon
^{\prime
})b+k(1-\lambda\varepsilon^{\prime})+t ) \,\mathrm{d} t}
\\
& \geq&\frac{\theta P ( \exists j\dvt X^{T}v_{j}>a((a_{j}+\lambda
\varepsilon^{\prime})b+k(1+\lambda\varepsilon^{\prime})) )
}{\int_{0}^{\infty}P ( \exists j\dvt X^{T}v_{j}>(\min_{j}a_{j}-\lambda
\varepsilon^{\prime})b+k(1-\lambda\varepsilon^{\prime})+t ) \,\mathrm{d} t}.
\end{eqnarray*}
We introduce the change of variables $t=s[(\min_{j}a_{j}-\lambda
\varepsilon^{\prime})b+k(1-\lambda\varepsilon^{\prime})]$ for the
integral in
the denominator and obtain that%
\begin{eqnarray*}
&& \int_{0}^{\infty}P \Bigl( \exists j\dvt
X^{T}v_{j}>\Bigl(\min_{j}a_{j}%
-\lambda\varepsilon^{\prime}\Bigr)b+k\bigl(1-\lambda\varepsilon^{\prime
}
\bigr)+t \Bigr) \,\mathrm{d} t
\\
&&\quad =\Bigl[\Bigl(\min_{j}a_{j}-\lambda
\varepsilon^{\prime}\Bigr)b+k\bigl(1-\lambda \varepsilon ^{\prime}
\bigr)\Bigr]\\
&&\qquad{}\times\int_{0}^{\infty}P \Bigl( \exists j\dvt
X^{T}v_{j}>\Bigl[\Bigl(\min_{j}
a_{j}-\lambda\varepsilon^{\prime}\Bigr)b+k\bigl(1-
\lambda\varepsilon^{\prime
}\bigr)\Bigr](1+s) \Bigr) \,\mathrm{d} s.
\end{eqnarray*}
Notice that for $b$ sufficiently large there exists $c^{\prime\prime
}$ so that%
\begin{eqnarray*}
&& \int_{0}^{\infty}P \Bigl( \exists j\dvt
X^{T}v_{j}>\Bigl[\Bigl(\min_{j}a_{j}%
-\lambda\varepsilon^{\prime}\Bigr)b+k\bigl(1-\lambda\varepsilon^{\prime
}
\bigr)\Bigr](1+s) \Bigr) \,\mathrm{d} s
\\
&&\quad \leq c^{\prime\prime}\theta P \bigl( \exists j\dvt X^{T}v_{j}>a
\bigl(\bigl(a_{j}%
+\lambda\varepsilon^{\prime}
\bigr)b+k\bigl(1+\lambda\varepsilon^{\prime
}\bigr)\bigr) \bigr) .
\end{eqnarray*}
Therefore, we have%
\[
p_{b}(s)\geq\frac{1}{c^{\prime\prime} ( (\min_{j}a_{j}-\lambda
\varepsilon)b+k(1-\lambda\varepsilon) ) }.
\]
The possibility of selecting $\delta^{\prime}>0$ to satisfy \eqref{Claimaa}
follows from the previous inequality.

Consequently, having (\ref{Claimaa}) in hand, if $t< c^{\prime
}\gamma$%
\begin{eqnarray*}
E_{0} \biggl( \prod_{k\leq tb}
\bigl(1-p_{b}(\hat{S}_{k})\bigr)I\bigl(\Vert \hat
{S}_{k}-\eta j\Vert _{2}<\varepsilon\max\{j,b\}\bigr)
\biggr) \leq\exp \biggl( -\sum_{1\leq j\leq tb
}
\frac{\delta^{\prime}}{j+b} \biggr) .
\end{eqnarray*}
We then conclude that if
%
\begin{equation}
\label{bbd}\mathop{\underline{\lim}}_{b\rightarrow\infty}P_{0} (
N_{b}%
<\infty ) \geq\mathop{\underline{\lim}}_{b\rightarrow\infty
}P_{0}
\bigl( N_{b}< c^{\prime}\gamma b \bigr) \geq1-\frac{1}{ ( c^{\prime}\gamma
+1 )
^{\delta^{\prime}}}
+\mathrm{o} (1)
\end{equation}
and this implies the statement of the result.
\end{pf}

\begin{lemma}
\label{LemTeqN2} For any $\varepsilon>0$ we can select $a,\gamma>0$
such that%
\[
\mathop{\underline{\lim}}_{b\rightarrow\infty}P_{0} ( T_{bA}=N_{b}
) \geq1-\varepsilon.
\]
\end{lemma}

\begin{pf}
Define%
\[
\tau ( -\delta_{2}b ) =\inf\bigl\{n\geq0\dvt r_{b}(S_{n})
\geq -\delta _{2}b\bigr\}.
\]
Let $c^{\prime}$ be as chosen in the proof of Lemma \ref{LemTeqN1}.
We then
have
\begin{eqnarray*}
&& P_{0}({T}_{bA}=N_{b},N_{b}<\infty)
\\
&&\quad =\sum_{k\leq c^{\prime}\gamma b}P_{0}\Bigl(N_{b}=k,
\max_{j=1}^{2m+d}\bigl[ {S}%
_{k}^{T}v_{j}-a_{i}b\bigr]
\geq0,T_{bA}>k-1\Bigr)
\\
&&\quad \geq\sum_{k\leq c^{\prime}\gamma b}P_{0}
\Bigl(N_{b}=k,\max_{j=1}^{2m+d}\bigl[\hat
{S}_{k}^{T}v_{j}-a_{i}b\bigr]\geq0,
\tau ( -\delta_{2}b ) >k-1\Bigr)
\\
&&\quad \geq(1-\varepsilon)\sum_{k\leq c^{\prime}\gamma
b}P_{0}
\bigl(N_{b}=k,\tau ( -\delta_{2}b ) >k-1\bigr).
\end{eqnarray*}
In the last inequality, we have used Lemma \ref{LemJump}. Further,%
\begin{eqnarray*}
\sum_{k\leq c^{\prime}\gamma b}P_{0}\bigl(N_{b}=k,
\tau ( -\delta _{2}b ) >k-1\bigr) & \geq& P_{0}
\bigl(N_{b}\leq c^{\prime}\gamma b,\tau ( -\delta
_{2}b ) >N_{b}-1\bigr)
\\
& \geq& P_{0}\bigl(N_{b}\leq c^{\prime}\gamma b,\tau
( -\delta _{2}b ) =\infty\bigr)
\\
& \geq& P \bigl( N_{b}\leq c^{\prime}\gamma b \bigr)
-P_{0}\bigl(\tau ( -\delta_{2}b ) <\infty\bigr).
\end{eqnarray*}
Since $P_{0}(\tau ( -\delta_{2}b ) <\infty)\rightarrow0$ as
$b\rightarrow\infty$, for $\gamma$ sufficiently large, we conclude
\begin{eqnarray*}
\mathop{\underline{\lim}}_{b\rightarrow\infty}P_{0} ( {T}_{bA}=N_{b}
) & \geq&\mathop{\underline{\lim}}_{b\rightarrow\infty
}P_{0}({T}_{bA}=N_{b},N_{b}<
\infty)
\\
& \geq&(1-\varepsilon)\mathop{\underline{\lim}}_{b\rightarrow\infty}P \bigl( N_{b}
\leq c^{\prime}\gamma b \bigr) \geq(1-\varepsilon)^{2},
\end{eqnarray*}
where the last equality follows from \eqref{bbd} in the proof of Lemma
\ref{LemTeqN1}. We conclude the proof.
\end{pf}


Proposition \ref{PropCLTAux} will follow as a consequence of the
next result.\vspace*{-1pt}

\begin{proposition}
\label{PropCoupling} By possibly enlarging the probability space, we
have the
following three coupling results.

\begin{longlist}[(iii)]
\item[(i)] Let $\kappa_{\mathbf{a}}$ be as defined in (\ref{DefKa}). There
exists a family of sets $(B_{b}\dvt b>0)$ such that $P ( B_{b} )
\rightarrow1$ as $b\nearrow\infty$ and with the property that if
$t\leq
\frac{\gamma}{2d}$
\[
P ( N_{b}>tb\vert S ) =P ( Z_{a,\theta
}>t ) \bigl( 1+\mathrm{o}  ( 1 )
\bigr)
\]
as $b\rightarrow\infty$ uniformly over $S\in B_{b}$, where
\[
P ( Z_{a,\theta}>t ) =\exp \biggl( -\theta\int_{0}^{t}
\frac
{\kappa_{a\mathbf{a}}(as)}{\int_{s}^{\infty}\kappa_{\mathbf
{a}}(u)\,\mathrm{d} u}%
\,\mathrm{d} s \biggr) .
\]

\item[(ii)] We can embed the random walk $S=\{S_{n}\dvt n\geq0\}$ and a uniform
random variable $U$ on $(0,1)$ independent of $S$ in a probability
space such
that $N_{b}\triangleq N_{b} ( S,U ) $ (a function of $S$ and $U$)
and $Z_{a,\theta}\triangleq Z_{a,\theta} ( U ) $ (a
function of
$U$) for all $S\in B_{b}$ such that
%
\begin{equation}
\biggl\llvert \frac{N_{b} ( S,U ) }{b}-Z_{a,\theta} ( U ) \biggr\rrvert
\rightarrow0 \label{Ap2}%
\end{equation}
as $b\rightarrow\infty$ for almost every $U\leq P(Z_{a,\theta}\leq
\frac
{\gamma}{2d})$. Furthermore, one can construct a $d$-dimensional Brownian
motion $B(t)$ so that%
%
\begin{equation}
S_{ \lfloor t \rfloor}=t\eta+CB ( t ) +e ( t ), \label{Sq}%
\end{equation}
where $CC^{T}=\operatorname{Var} ( X ) $ is the covariance matrix of an increment
$X$ and $e ( \cdot ) $ is a (random) function such that
\[
\frac{|e ( xt ) |}{t^{1/2}}\rightarrow0
\]
with probability one, uniformly on compact sets on $x\geq0$ as
$t\rightarrow
\infty$.

\item[(iii)] Finally, we can also embed a family of random variables $\{
\hat
{X} ( a,b,s ) \}$, independent of $S$ and $U$, distributed
as $X$
conditioned on $s+X\in A_{b,a} ( s ) $, coupled with a random
variable $Y ( z ) $ (also independent of $S$ and $U$) so
that for
each Borel set $B$
%
\[
P \bigl( Y ( z ) \in B \bigr) =\frac{\mu ( B\cap
\{y\dvt \max_{j=1}^{2m+d}(y^{T}v_{j}-a_{i})\geq z\} ) }{\mu (
\{y\dvt \max_{j=1}^{2m+d}(y^{T}v_{j}-a_{i})\geq z\} ) },
\]
and with the property%
\[
\lim_{a\rightarrow1,b\rightarrow\infty}\biggl\vert\frac{\hat{X} (
a,b,\eta ( z+\xi_{b} )  ) }{b}-Y ( z ) \biggr\vert \rightarrow0
\]
with probability 1 as long as $\xi_{b}\rightarrow0$.
\end{longlist}
\end{proposition}

\begin{pf}
We start with the proof of (i). Note that%
\[
P_{0} ( N_{b}>tb\vert S ) =\prod
_{0\leq
k\leq \lfloor tb \rfloor-1}\bigl(1-p_{b} ( S_{j} ) \bigr),
\]
where, by convention, a product indexed by an empty subset is equal to unity.
Now, let $\delta_{b}=1/\log b$, $\gamma_{k,\delta_{b}}=\max (
1/\delta_{b}^{2},\delta_{b}k ) $ and
\[
B_{b}= \bigl\{ S\dvt \Vert S_{k}-k\eta
\Vert _{2}\leq\gamma_{k,\delta_{b}} \bigr\}
\]
for $k\geq1$ for all $k\leq tb$. It follows easily that
$P(B_{b})\rightarrow1$
as $b\rightarrow\infty$. Recall that if $r_{b} ( s ) \leq
-\delta_{2}b$ and $s^{T}\eta\leq\gamma b$%
\[
p_{b}(s)  =\frac{\theta P ( \exists j\dvt X^{T}v_{j}>a(a_{j}b-s^{T}%
v_{j}) ) }{\int_{0}^{\infty}P ( \exists j\dvt X^{T}v_{j}>a_{j}%
b-s^{T}v_{j}+t ) \,\mathrm{d} t}.
\]
We will find upper and lower bounds on the numerator and denominator on the
set $B_{b}$ so that we can use regular variation properties to our advantage.
Suppose that $S_{k} = s$. Define $\lambda=\max_{j}\Vert v_{j}%
\Vert _{2}$ and observe that if $\Vert
s-k\eta
\Vert _{2}\leq\gamma_{k,\delta_{b}}$ then for
all $1\leq
j\leq2m+d$ and all $k\geq0$ we have that%
\[
\bigl\llvert s^{T}v_{j}-k\eta^{T}v_{j}
\bigr\rrvert =\bigl\llvert s^{T}%
v_{j}+k\bigr
\rrvert \leq\lambda\gamma_{k,\delta_{b}}.
\]
In addition,
\[
s^{T}\eta\leq kd+d^{1/2} \gamma_{k,\delta_{b}}\leq\gamma b
\]
if $k\leq tb$ and $t\leq\frac{\gamma}{2d}$, for all $b$ large enough.
Therefore, if $X^{T}v_{j}>a(a_{j}b-s^{T}v_{j})$ and $\Vert
s-k\eta\Vert _{2}\leq\gamma_{k,\delta_{b}}$,%
\[
X^{T}v_{j}>a\bigl(a_{j}b-s^{T}v_{j}
\bigr)\Rightarrow X^{T}v_{j}>aa_{j}b+ak-\lambda
\gamma_{k,\delta_{b}}.
\]
Consequently, if $\Vert  s-k\eta\Vert _{2}%
\leq\gamma_{k,\delta_{b}}$%
\begin{eqnarray*}
&&P \bigl( \exists j\dvt X^{T}v_{j}>a\bigl(a_{j}b-s^{T}v_{j}
\bigr) \bigr) \\
&&\quad\leq P \bigl( \exists j\dvt X^{T}v_{j}>aa_{j}b+ak-a
\lambda\gamma_{k,\delta_{b}} \bigr) .
\end{eqnarray*}
Similarly,%
\begin{eqnarray*}
&&P \bigl( \exists j\dvt X^{T}v_{j}>a\bigl(a_{j}b-s^{T}v_{j}
\bigr) \bigr) \\
&&\quad\geq P \bigl( \exists j\dvt X^{T}v_{j}>a_{j}b+ak+a
\lambda\gamma_{k,\delta_{b} } \bigr) .
\end{eqnarray*}
Following analogous steps, we obtain that for $\llVert  s-k\eta
\rrVert
_{2}\leq\gamma_{k,\delta_{b}}$%
\[
v_{b}(s)\leq\int_{0}^{\infty}P \bigl(
\exists j\dvt X^{T}v_{j}>a_{j}%
b+k+t-
\lambda\gamma_{k,\delta_{b}} \bigr) \,\mathrm{d} t
\]
and%
\[
v_{b}(s)\geq\int_{0}^{\infty}P \bigl(
\exists j\dvt X^{T}v_{j}>a_{j}%
b+k+t+
\lambda\gamma_{k,\delta_{b}} \bigr) \,\mathrm{d} t.
\]
Our previous upper and lower bounds indicate that on $\llVert  S_{k}%
-k\eta\rrVert _{2}\leq\gamma_{k,\delta_{b}}$,%
\[
v_{b}(S_{k})=\bigl(1+\mathrm{o}  ( 1 ) \bigr)\int
_{0}^{\infty}P \bigl( \exists j\dvt
X^{T}v_{j}>a_{j}b+k+t \bigr) \,\mathrm{d} t
\]
on the set $B_{b}$, uniformly as $b\rightarrow\infty$ and then%
\begin{eqnarray*}
p_{b}(S_{k}) & =&\theta\bigl(1+\mathrm{o}  ( 1 ) \bigr)
\frac{P ( \exists
j\dvt X^{T}v_{j}>aa_{j}b+ak ) }{\int_{0}^{\infty}P ( \exists
j\dvt X^{T}v_{j}>a_{j}b+k+t ) \,\mathrm{d} t}
\\
& =&\theta\bigl(1+\mathrm{o}  ( 1 ) \bigr)\frac{P ( \exists j\dvt X^{T}v_{j}%
>aa_{j}b+ak ) }{\int_{k}^{\infty}P ( \exists j\dvt X^{T}v_{j}%
>a_{j}b+u ) \,\mathrm{d} u}
\\
& =&\theta\bigl(1+\mathrm{o}  ( 1 ) \bigr)\frac{P ( \exists j\dvt X^{T}v_{j}%
>aa_{j}b+ak ) }{b\int_{k/b}^{\infty}P ( \exists
j\dvt X^{T}v_{j}%
>a_{j}b+tb ) \,\mathrm{d} t}.
\end{eqnarray*}
By definition of regular variation there exists a slowly varying function
$L ( \cdot ) $ such that%
\[
P \bigl( \exists j\dvt X^{T}v_{j}-aa_{j}b>z
\bigr) =L ( b ) b^{-\alpha}\kappa_{a\mathbf{a}} ( z/b ) .
\]
Therefore, by Karamata's theorem (see \cite{RES87}) we then have that
if $S\in
B_{b}$ and $k\leq tb$ with $t\leq\gamma/2d$%
\[
p_{b}(S_{k})=\theta\bigl(1+\mathrm{o}  ( 1 ) \bigr)
\frac{\kappa_{a\mathbf
{a}} (
ak/b ) }{b\int_{k/b}^{\infty}\kappa_{\mathbf{a}} (
t ) \,\mathrm{d} t}%
\]
as $b\rightarrow\infty$ uniformly over $S$ in $B_{b}$. Therefore,
\begin{eqnarray*}
P_{0} ( N_{b}>tb\vert S )  =\prod
_{0\leq
j\leq \lfloor tb \rfloor-1}\bigl(1-p_{b} ( S_{j} ) \bigr)
 =\exp \Biggl( - \bigl( \theta+\mathrm{o}  ( 1 ) \bigr) \sum
_{j=0}^{ \lfloor tb \rfloor-1}\frac{\kappa_{a\mathbf
{a}} (
ak/b ) }{b\int_{k/b}^{\infty}\kappa_{\mathbf{a}} (
t )
\,\mathrm{d} t} \Biggr)
\end{eqnarray*}
uniformly as $b\nearrow\infty$ over $S\in B_{b}$ and
\[
\sum_{j=0}^{ \lfloor tb \rfloor-1}\frac{\kappa_{a\mathbf
{a}} (
ak/b ) }{b\int_{k/b}^{\infty}\kappa_{\mathbf{a}} (
t )
\,\mathrm{d} t}
\rightarrow\int_{0}^{t}\frac{\kappa_{a\mathbf{a}} (
as ) }%
{\int_{s}^{\infty}\kappa_{\mathbf{a}}(u)\,\mathrm{d} u}\,\mathrm{d} s,
\]
which implies that as long as $t\leq\frac{\gamma}{2d}$%
\[
P_{0} ( N_{b}>tb\vert S ) \rightarrow P (
Z_{a,\theta}>t ) =\exp \biggl( -\theta\int_{0}^{t}
\frac
{\kappa
_{a\mathbf{a}} ( as ) }{\int_{s}^{\infty}\kappa_{\mathbf
{a}}%
(u)\,\mathrm{d} u}\,\mathrm{d} s \biggr) .
\]
Part (ii) is straightforward from (i) by Skorokhod embedding for random
walks as
follows. First, construct a probability space in which the random walk
$S$ is
strongly approximated by a Brownian motion as in (\ref{Sq}) and
include a
uniform random variable $U$ independent of $S$. Construct $N_{b}/b$ applying
the generalized inverse cumulative distribution function of $N_{b}/b$ given
$S$ to $U$. Then, apply the same uniform $U$ to generate $Z_{a,\theta
}$ by
inversion. Because of our estimates in part (i) we have that (\ref
{Ap2}) holds
as long as $U\leq P(Z_{a,\theta}\leq\frac{\gamma}{2d})$. Part (iii) follows
using a similar argument.
\end{pf}

Now we are ready to provide the proof of Proposition \ref{PropCLTAux}.

\begin{pf*}{Proof of Proposition \ref{PropCLTAux}}
According to Proposition
\ref{PropFinieH}, the distribution of $\hat{S}$ approximates the conditional
random walk up to time $T_{bA}-1$ in total variation. By virtue of Lemmas
\ref{LemTeqN1} and \ref{LemTeqN2}, we can replace
$(T_{bA},X_{T_{bA}})$ by
$(N_{b},\hat{X}_{N_{b}})$. Then, it suffices to show weak convergence of
\[
\biggl( \frac{N_{b}}{b},\frac{S_{uN_{b}}-uN_{b}\eta}{\sqrt {N_{b}}},\frac
{\hat{X}_{N_{b}}}{b} \biggr) ,
\]
given that $\theta$ and $a$ can be chosen arbitrarily close to 1. By
using the
embedding in Proposition \ref{PropCoupling} we consider $S\in B_{b}$
so that
on $U\leq P(Z_{a,\theta}\leq\gamma/2)$%
\begin{eqnarray*}
&& \biggl( \frac{N_{b}}{b},\frac{S_{uN_{b}}-uN_{b}\eta}{\sqrt{N_{b}}} 
,
\frac{\hat X_{N_{b}}}{b} \biggr)
\\
&&\quad = \biggl( Z_{a,\theta}+\xi_{b},\frac{CB ( ubZ_{a,\theta
}+ub\xi
_{b} ) +e ( uN_{b} ) }{\sqrt{bZ_{a,\theta}+\xi_{b}}}
,Y_{a} ( Z_{a,\theta}+\xi_{b} ) +
\chi_{b} \biggr) ,
\end{eqnarray*}
where $\xi_{b},\chi_{b}\rightarrow0$ as $b\rightarrow\infty$ and
$\sup_{u\in\lbrack0,1]}|e(uN_{b})/\sqrt{N_{b}}|\rightarrow0$ as $N_{b}%
\rightarrow\infty$. Also, note that as $b\rightarrow\infty$ we
choose $a$ and
$\theta$ close $1$ and $\gamma$ sufficiently large. From Lemmas \ref
{LemTeqN1}
and \ref{LemTeqN2}, we must verify that for each $z>0$
\[
\sup_{0\leq u\leq1}\biggl\llvert \frac{B ( ubz+ub\xi_{b} )
-B (
ubz ) }{b^{1/2}}\biggr\rrvert
\longrightarrow0.
\]
Given $z>0$ select $b$ large enough so that $\xi_{b}\leq\varepsilon$
for each
$S\in B_{b}$. Then, it suffices to bound the quantity%
\[
\sup_{u,s\in(0,1)\dvt \llvert  u-s\rrvert \leq\varepsilon}\biggl\llvert \frac{B ( ubz ) -B ( sbz ) }{b^{1/2}}\biggr\rrvert .
\]
However, by the invariance principle, the previous quantity equals in
distribution to%
\[
z^{1/2}\sup_{u,s\in(0,1)\dvt \llvert  u-s\rrvert \leq\varepsilon
}\bigl\llvert B ( u ) -B ( s ) \bigr
\rrvert ,
\]
which is precisely the modulus of continuity of Brownian motion evaluated
$\varepsilon$, which (by continuity of Brownian motion) goes to zero almost
surely as $\varepsilon\rightarrow0$.
From Proposition \ref{PropCoupling}, we have that $Z_{a,\theta
}\rightarrow Z$
as $a,\theta\rightarrow1$, where $Z$ is defined as in the statement
of the
proposition. Since $\gamma$ can be chosen arbitrarily large as
$b\rightarrow
\infty$, we complete the proof.
\end{pf*}

Finally, we give the proof of Theorem \ref{ThmMain2}.

\begin{pf*}{Proof of Theorem \ref{ThmMain2}}It suffices to exhibit a coupling
under which%
\[
\bigl( Z,CB ( uZ ) ,Y ( Z ) \bigr) - \bigl( Z^{\ast
},CB \bigl(
uZ^{\ast} \bigr) ,Y^{\ast} \bigl( Z^{\ast} \bigr) \bigr)
\longrightarrow0
\]
almost surely as $\beta,\gamma\rightarrow\infty$ and $\delta
\rightarrow0$, but
this is immediate from continuity of the Brownian motion and of the radial
component of the measure $\mu ( \cdot ) $ (Lemmas \ref{LemPropmu}
and~\ref{LemRatv}).
\end{pf*}

\begin{appendix}
\section*{Appendix: Some properties of regularly varying
distributions}\label{SectionAppendixProperties}

In this Appendix, we summarize some important properties of the regularly
varying distribution of $X$, which satisfies the assumptions stated in Section
\ref{SecISPS}. We are mostly concerned with some continuity
properties of
the limiting measure $\mu ( \cdot ) $ (defined in equation
(\ref{RVDist})).

The measure $\mu ( \cdot ) $ can be
represented as a \textit{product measure} corresponding to the angular
component and the radial component \cite{RES06}. The angular component
\[
\Phi( \cdot) =\frac{\mu(\{x\dvt \Vert x\Vert_2 >1 ,x/\Vert x\Vert_2
\in\cdot\})}{\mu(\{x\dvt \Vert x\Vert_2 >1 \})}
\]
corresponds to a probability measure on
the $(d-1)$-dimensional sphere in $\mathbb{R }^{d}$. The radial
component $\vartheta ( \mathrm{d}r ) $ is a measure that is
absolutely continuous with respect to the Lebesgue measure. Moreover,
$\vartheta ( \mathrm{d}r ) =cr^{-\alpha-1}\,\mathrm{d} r$ for some constant
$c>0$. Then, the measure $\mu$ can be written as the product of
$\vartheta$ and $\Phi$.
We then obtain the following lemma.

\begin{lemma}
\label{LemPropmu}Let $\Gamma=\{y\dvt \eta^{T}y>\gamma\}$ and $R_{2}=\{
y\dvt \max_{i=1}^{2m+d} ( y^{T}v_{i}-a_{i} ) \leq-\delta_{2}\}$ for
$\gamma>0$ large but fixed and $\delta_{2}>0$ small. Let $K=\Gamma
^{c}\cup
R_{2}$ and define $\kappa ( \cdot ) $ for each $t>0$,
$z\in K$ and
$a^{\prime}=(a_{1}^{\prime},\ldots,a_{2m+d}^{\prime})$ in a small
neighborhood of
$a= ( a_{1},\ldots,a_{2m+d} ) $ via%
\[
\kappa_{\mathbf{a}^{\prime}} ( t,z ) =\mu \Bigl( \Bigl\{ y\dvt \max
_{j=1}^{2m+d} \bigl( y^{T}v_{j}+z^{T}v_{j}-a_{j}^{\prime}
\bigr) >t\Bigr\} \Bigr) .
\]
Then,%
\renewcommand{\theequation}{\arabic{equation}}
\setcounter{equation}{66}
\begin{equation}
\kappa_{\mathbf{a}^{\prime}} ( t,z ) =\int_{\mathcal
{S}_{d}}%
\frac{c}{\alpha}\max_{j=1}^{2m+d} \biggl(
\frac{\cos ( \theta
,v_{j} ) ^{+}}{-z^{T}v+a_{j}^{\prime}+t} \biggr) ^{\alpha}\Phi ( \mathrm{d}\theta ) .\label{RepLemMu}%
\end{equation}

\end{lemma}

\begin{pf}The result follows immediately from the representation of $\mu
(
\cdot ) $ in polar coordinates. We shall sketch the details. For each
$j$, let $\mathcal{H}_{j} ( a_{j}^{\prime},t ) =\{y\dvt y^{T}%
v_{j}+z^{T}v_{j}-a_{j}^{\prime}>t\}$ and note that in polar
coordinates, we can
represent $\mathcal{H}_{j} ( a_{j}^{\prime},t ) $ as%
\[
\bigl\{ ( \theta,r ) \dvt r>0\mbox{ and }r\Vert
v_{j}\Vert _{2}\cos ( \theta,v_{j}
) \geq -z^{T}v+a_{j}^{\prime}+t\bigr\}.
\]
Note that $-z^{T}v+a_{j}^{\prime}+t>0$ for all $t\geq0$, $j\in\{
1,\ldots,2m+d\}$
and $z\in K$. Therefore,
\begin{eqnarray*}
&& \mu \Bigl( \Bigl\{y\dvt \max_{j=1}^{2m+d} \bigl(
y^{T}v_{j}+z^{T}v_{j}-a_{j}%
^{\prime} \bigr) >t\Bigr\} \Bigr)
\\
&&\quad =\int_{\mathcal{S}_{d}}\int_{0}^{\infty}cr^{-\alpha-1}I
\Bigl(r\geq \min_{j}\bigl\{\bigl(-z^{T}v_{j}+a_{j}^{\prime}+t
\bigr)/\cos ( \theta,v_{j} ) ^{+}\bigr\}\Bigr)\,\mathrm{d} r\Phi ( \mathrm{d}
\theta )
\\
&&\quad =\int_{\mathcal{S}_{d}}\frac{c}{\alpha}\max_{j=1}^{2m+d}
\biggl( \frac
{\cos ( \theta,v_{j} ) ^{+}}{-z^{T}v_{j}+a_{j}^{\prime
}+t} \biggr) ^{\alpha}\Phi ( \mathrm{d}\theta ) ,
\end{eqnarray*}
and the result follows.\vadjust{\goodbreak}
\end{pf}

\begin{lemma}
\label{LemRatv}Assume that $K=\Gamma^{c}\cup R_{2}$ is defined as
in the
previous lemma and note that $K$ is a non-empty compact set. Suppose that
$z\in K$, define $s=zb$ and write%
\begin{eqnarray*}
v_{b} ( zb ) & =&\int_{0}^{\infty}P \Bigl(
\max_{i=1}%
^{2m+d} \bigl(
v_{i}^{T}X+z^{T}v_{i}b-a_{i}b
\bigr) >t \Bigr) \,\mathrm{d} t,
\\
\kappa_{\mathbf{a}} ( t,z ) & =&\mu \Bigl( \Bigl\{y\dvt \max
_{j=1}%
^{2m+d} \bigl( y^{T}v_{j}+z^{T}v_{j}-a_{j}
\bigr) >t\Bigr\} \Bigr) .
\end{eqnarray*}
Then,
\[
\lim_{b\rightarrow\infty}\sup_{z\in K}\biggl\llvert
\frac{v_{b}
( zb )
}{bP ( \Vert  X\Vert
_{2}>b )
\int_{0}^{\infty}\kappa_{\mathbf{a}} ( t,z ) \,\mathrm{d} t}-1\biggr\rrvert \longrightarrow0
\]
as $b\rightarrow\infty$.
\end{lemma}

\begin{pf}We first write%
\[
v_{b} ( zb ) =b\int_{0}^{\infty}P \Bigl(
\max_{i=1}^{2m+d} \bigl( v_{i}^{T}X+bz^{T}v_{i}-a_{i}b
\bigr) >ub \Bigr) \,\mathrm{d} u,
\]
and define%
\[
\mathcal{L}_{b} ( u,z ) =P \Bigl( \max_{i=1}^{2m+d}
\bigl( v_{i}^{T}X+bz^{T}v_{i}-a_{i}b
\bigr) >ub \Bigr) .
\]
Set $\varepsilon\in ( 0,\delta_{2} ) $ arbitrarily small
but fixed
and let $M=\max_{i\leq2m+d}\Vert v_{i}\Vert _{2} <0$. Define $\varepsilon
^{\prime
}=\varepsilon/M$ and consider an open cover of the set $K$ by balls with
radius $\varepsilon^{\prime}$ centered at points $\mathcal
{C}_{\varepsilon
^{\prime}}=\{w_{1},\ldots,w_{m^{\prime}}\}\subset K$. We then have that
for every
$z\in K$ there exists $w_{k}\triangleq w_{k} ( z ) \in
\mathcal{C}_{\varepsilon^{\prime}}$ such that $\Vert
z-w_{k}\Vert _{2}\leq\varepsilon^{\prime}$. Note
that for
each $z\in K$
\[
\mathcal{L}_{b} ( u+\varepsilon,w_{k} ) \leq\mathcal
{L}_{b} ( u,z ) \leq\mathcal{L}_{b} ( u-
\varepsilon,w_{k} ) .
\]
Consequently,%
\[
\int_{0}^{\infty}\mathcal{L}_{b} ( u+
\varepsilon,w_{k} ) \,\mathrm{d} u\leq v_{b} ( zb ) \leq\int
_{0}^{\infty}\mathcal{L}_{b} ( u-
\varepsilon,w_{k} ) \,\mathrm{d} u.
\]
Now we claim that%
\renewcommand{\theequation}{\arabic{equation}}
\setcounter{equation}{67}
\begin{equation}
\lim_{b\rightarrow\infty}\frac{\int_{0}^{\infty}\mathcal
{L}_{b} (
u+\varepsilon,w_{k} ) \,\mathrm{d} u}{bP ( \Vert
X\Vert _{2}>b ) \int_{0}^{\infty}\kappa_{\mathbf
{a}} (
t+\varepsilon,w_{k} ) \,\mathrm{d} t}=1. \label{Claimlim}%
\end{equation}
The previous limit follows from dominated convergence as follows.
First, we
have that%
\[
\frac{\mathcal{L}_{b} ( u+\varepsilon,w_{k} ) }{bP (
\Vert  X\Vert _{2}>b )
}\longrightarrow \kappa_{\mathbf{a}} ( u+
\varepsilon,w_{k} )
\]
for every $u$ fixed by the definition regular variation. Then, if $u\in(0,C)$
for any $C>0$ we conclude that%
\[
\frac{\mathcal{L}_{b} ( u+\varepsilon,w_{k} ) }{bP (
\Vert  X\Vert _{2}>b ) }\leq \frac{\mathcal{L}_{b} ( 0,w_{k} ) }{bP ( \Vert
X\Vert _{2}>b ) }=\mathrm{O}  ( 1 )
\]
as $b\rightarrow\infty$ and therefore by the bounded convergence
theorem, we
conclude that%
\[
\lim_{b\rightarrow\infty}\frac{\int_{0}^{C}\mathcal{L}_{b} (
u+\varepsilon,w_{k} ) \,\mathrm{d} u}{bP ( \Vert
X\Vert _{2}>b ) }=\int_{0}^{C}
\kappa_{\mathbf{a}} ( t+\varepsilon,w_{k} ) \,\mathrm{d} t.
\]
On the set $u\geq C$ we have that if $c_{i,k}=a_{i}+\llvert
w_{k}^{T}%
v_{i}\rrvert $, then%
\[
\mathcal{L}_{b} ( u+\varepsilon,w_{k} ) \leq\sum
_{i=1}%
^{2m+d}P \bigl( \Vert
X\Vert _{2}\Vert v_{i}\Vert _{2}> ( u+\varepsilon ) b-c_{i,k}b \bigr) .
\]
By Karamata's theorem for one dimensional regularly varying random variables,
it follows that if $\alpha>1$, then for $C\geq\max_{i,k}c_{i,k}$ the
functions
\[
\frac{P ( \Vert  X\Vert
_{2}\Vert  v_{i}\Vert _{2}> ( \cdot
+\varepsilon )
b-c_{i,k}b ) }{bP ( \Vert  X\Vert
_{2}>b ) }%
\]
are uniformly integrable with respect to the Lebesgue measure on
$[C,\infty)$
and therefore we conclude that%
\begin{eqnarray*}
&&\lim_{b\rightarrow\infty}\frac{\int_{C}^{\infty}\mathcal
{L}_{b} (
u+\varepsilon,w_{k} ) \,\mathrm{d} u}{bP ( \Vert
X\Vert _{2}>b ) }\\
&&\quad=\int_{C}^{\infty}
\kappa_{\mathbf
{a}} ( t+\varepsilon,w_{k} ) \,\mathrm{d} t
\end{eqnarray*}
and therefore the limit (\ref{Claimlim}) holds. Thus, we have that%
\begin{eqnarray*}
&& \sup_{\{z\in\Vert  z-w_{k}\Vert _{2}
\leq\varepsilon^{\prime}\}}\biggl\llvert \frac{\int_{0}^{\infty
}\mathcal{L}%
_{b} ( u+\varepsilon,w_{k} ) \,\mathrm{d} u}{bP ( \Vert
X\Vert _{2}>b ) \int_{0}^{\infty}\kappa
_{\mathbf{a}%
} ( t,z ) \,\mathrm{d} t}-1\biggr\rrvert
\\
&&\quad =\sup_{\{z\in\Vert  z-w_{k}\Vert _{2}
\leq\varepsilon^{\prime}\}}\biggl\llvert \frac{\int_{0}^{\infty
}\kappa
_{\mathbf{a}} ( t+\varepsilon,z ) \,\mathrm{d} t}{\int_{0}^{\infty}%
\kappa_{\mathbf{a}} ( t,z ) \,\mathrm{d} t}-1\biggr\rrvert
+\mathrm{o}  ( 1 )
\end{eqnarray*}
as $b\rightarrow\infty$. Observe from representation (\ref
{RepLemMu}) and
Assumption \ref{ass2} in Section \ref{SecISPS} we have that%
\[
\kappa_{\mathbf{a}} ( t,z ) >\mu \Bigl( \Bigl\{y\dvt \max
_{j=1}^{2m+d}%
y^{T}v_{j}>
\delta_{2}\Bigr\} \Bigr) >0.
\]
Moreover, it also follows as an easy application of the dominated convergence
theorem and our representation in (\ref{RepLemMu}) that%
\[
\lim_{\varepsilon\rightarrow0}\sup_{\{z\in\Vert
z-w_{k}%
\Vert _{2} \leq\varepsilon^{\prime}\}}\biggl\llvert
\frac
{\int_{0}^{\infty}\kappa_{\mathbf{a}} ( t+\varepsilon,z
) \,\mathrm{d} t}%
{\int_{0}^{\infty}
\kappa_{\mathbf{a}} ( t,z ) \,\mathrm{d} t}-1\biggr\rrvert =0.
\]
We then conclude that%
\begin{eqnarray*}
&& \mathop{\overline{\lim}}_{b\rightarrow\infty}\sup_{z\in K}\biggl\llvert
\frac{\int_{0}^{\infty}\mathcal{L}_{b} ( u,z ) \,\mathrm{d} u}{bP ( \Vert  X\Vert _{2}>b ) \int_{0}^{\infty
}%
\kappa_{\mathbf{a}} ( t,z ) \,\mathrm{d} t}-1\biggr\rrvert
\\
&&\quad \leq\max_{k}\mathop{\overline{\lim}}_{b\rightarrow\infty}\sup
_{\{z\in
\Vert  z-w_{k}\Vert _{2} \leq\varepsilon
^{\prime}%
\}}\biggl\llvert \frac{\int_{0}^{\infty}\mathcal{L}_{b} (
u+\varepsilon
,w_{k} ) \,\mathrm{d} u}{bP ( \Vert  X\Vert
_{2}>b ) \int_{0}^{\infty}\kappa_{\mathbf{a}} (
t,z )
\,\mathrm{d} t}-1\biggr\rrvert
\\
& &\qquad{}+\max_{k}\mathop{\overline{\lim}}_{b\rightarrow\infty}\sup
_{\{z\in\Vert  z-w_{k}\Vert _{2} \leq\varepsilon
^{\prime}%
\}}\biggl\llvert \frac{\int_{0}^{\infty}\mathcal{L}_{b} (
u+\varepsilon
,w_{k} ) \,\mathrm{d} u-\int_{0}^{\infty}\mathcal{L}_{b} (
u-\varepsilon
,w_{k} ) \,\mathrm{d} u}{bP ( \Vert  X\Vert
_{2}>b ) \int_{0}^{\infty}\kappa_{\mathbf{a}} (
t,z )
\,\mathrm{d} t}\biggr\rrvert
\\
&&\quad \leq\max_{k}\biggl\llvert \frac{\int_{0}^{\infty}\kappa_{\mathbf
{a}} (
t+\varepsilon,w_{k} ) \,\mathrm{d} t}{\int_{0}^{\infty}\kappa_{\mathbf
{a}} (
t,z ) \,\mathrm{d} t}-1\biggr
\rrvert +\max_{k}\biggl\llvert \frac{\int_{0}^{\infty}%
\kappa_{\mathbf{a}} ( t+\varepsilon,w_{k} ) \,\mathrm{d} t-\int_{0}^{\infty
}\kappa_{\mathbf{a}} ( t-\varepsilon,w_{k} ) \,\mathrm{d} t}{\int_{0}^{\infty
}\kappa_{\mathbf{a}} ( t,z ) \,\mathrm{d} t}\biggr
\rrvert .
\end{eqnarray*}
Once again use the representation in Lemma \ref{LemPropmu} and the dominated
convergence theorem to conclude that the right-hand side of the previous
inequality can be made arbitrarily small as $\varepsilon\rightarrow
0$, thereby
concluding our result.
\end{pf}
\end{appendix}

\section*{Acknowledgements}
We are grateful to the Associate Editor and the referees for their helpful comments.
This research is supported in part by Institute of Education Sciences,
through Grant R305D100017, NSF DMS-08-06145, NSF CMMI-0846816, and NSF
CMMI-1069064.

%


\printhistory

\end{document}